\documentclass[12pts]{article}
\usepackage[utf8]{inputenc}
\usepackage{amsfonts}
\usepackage{amsmath}
\usepackage{amssymb}
\usepackage{fullpage}
\usepackage{hyperref}
\usepackage{geometry}
\usepackage{amsthm}
\usepackage{setspace}
\usepackage{multirow}
\usepackage{graphicx}
\usepackage{float}
\usepackage[affil-it]{authblk}
\newtheorem{theorem}{Theorem}[section]

\usepackage{xcolor}

\newtheorem{definition}{Definition}
\title{Model Identifiability for Bivariate Failure Time Data with Competing Risk: Non-parametric Cause-specific Hazards and Gamma Frailty}
\author[1]{Biswadeep Ghosh}
\author[1]{Anup Dewanji}
\author[2]{Sudipta Das}
\affil[1]{Applied Statistics Unit, Indian Statistical Institute}
\affil[2]{Ramakrishna Mission Vivekananda Educational and Research Institute}
\date{}

\begin{document}
	

\maketitle
{\bf Abstract.} In survival analysis, frailty variables are often used to model the association in multivariate survival data. Identifiability is an important issue while working with such multivariate survival data with or without competing risks. In this work, we consider bivariate survival data with competing risks and investigate identifiability results with non-parametric baseline cause-specific hazards and different types of Gamma frailty. Prior to that, we prove that, when both baseline cause-specific hazards and frailty distributions are non-paramteric, the model is not identifiable. We also construct a non-identifiable model when baseline cause-specific hazards are non-parametric but frailty distribution may be parametric. Thereafter, we consider four different Gamma frailty distributions and the corresponding models are shown to be identifiable under fairly general assumptions.

\section{Introduction}

Bivariate failure time data arises when, for example, an individual gives pair of observations (e.g., onset of hearing disability in the two ears), or two related individuals give the same end-point (onset of puberty in a study of twins). Modeling and analysis of bivariate failure time data have been well-studied in the recent literature. See Yashin et al.$(1995)$, Hougaard$(2000)$, 
Wienke$(2010)$, 
Iachine$(2004)$, among others, for example. 
We shall refer to the two components of the bivariate failure time as the two individuals in a pair. The dependence between the two failure times has been modeled in various ways. One way is to define bivariate cumulative hazard vector, as in Johnson and Kotz (1975) and Dabrowska (1988). Recently, there are two popular approaches to model the bivariate failure time distribution, one through the use of copula (See, for example, Prenen et. al, $(2017)$ and the other through the introduction of random frailty (See, for example,  Duchateau, $2008$ and Wienke, $2010$). Modeling through frailty has somehow become more common than the use of copula, possibly because of its appealing interpretation. \\

In this work, we consider bivariate failure time with competing risks. There have been some work on this using bivariate cumulative hazard vector (See Sankaran et al., 2006). Generalization of the frailty approach to multivariate failure time with competing risks seems to be very natural and has recently been studied by Gorfine and Hsu $(2011)$. We do not know of any copula based modeling approach of bivariate failure time with competing risks. This also does not seem to be that natural. \\

Identifiability is an important aspect of model building, specially in the context of frailty modeling, that has not been studied much. In particular, we shall show that the model is not identifiable when both the cause-specific hazards and the frailty distributions are unspecified. There are also examples of non-parametric cause-specific hazards and parametric frailty resulting in non-identifiable models. In this paper, we consider non-parametric cause-specific hazards and gamma frailty to investigate identifiability of the corresponding model for bivariate failure time with competing risks. In particular, we consider four different types of Gamma frailty depending on the association at individual and cause-specific level, namely, (1) shared Gamma frailty, (2) correlated Gamma frailty, (3) shared cause-specific Gamma frailty and (4) correlated cause-specific Gamma frailty, respectively. See Ghosh et al. (2024) for more details. \\

In Section 2, we develop the frailty model for bivariate failure time with competing risks in terms of the joint distribution function. This section also formally introduces the concept of identifiability and proves non-identifiability when both the cause-specific hazards and the frailty distributions are arbitrary. Section 3 presents an example of a non-identifiable model even when the frailty distribution is assumed to be parametric. Section 4 studies identifiability for shared Gamma frailty, while Section 5 does the same for correlated Gamma frailty. Sections 6 and 7 consider shared and correlated cause-specific Gamma frailty, respectively. Section 8 ends with some concluding remarks. \\

\vskip10pt

\section{Modeling}

\vskip10pt

A typical observation for bivariate survival data with competing risks for a pair of individuals is of the form $(T_1,T_2,J_1,J_2)$, where $T_k$ is the failure time of the $k$th individual of a pair with $J_k$ as the cause of failure, for $k=1,2$. Here, $T_1$ and $T_2$ are positive random variables and the random variable $J_k$ has the support $\{1,\cdots L_k\}$, the set of $L_k$ risks, only one of which is responsible for failure of the $k$th individual, for $k=1,2$.  
Bivariate survival data with competing risks can be modeled through the joint sub-distribution function $F_{j_1j_2}(t_1,t_2)$ defined as
$$F_{j_1j_2}(t_1,t_2) = P\big[T_1 \leq t_1, T_2 \leq t_2, J_1 = j_1, J_2 = j_2\big], $$
where $t_{k} > 0, j_{k}=1,\cdots,L_k, k = 1,2$. The joint survival function $S(t_1,t_2)$ is defined as
$$S(t_1,t_2) = P\big[T_1 > t_1,T_2 > t_2\big].$$ 
Similarly, the marginal sub-distribution function $F^{(k)}_{j}(t_k)$ is defined as  
$$F^{(k)}_{j}(t_k) = P\big[T_{k} \leq t_{k},J_{k} = j\big], \quad\mbox{ for }\quad t_{k} > 0,j=1,\cdots,L_k,\quad\mbox{ and }\quad k = 1,2.$$

We now build several frailty models for $F_{j_1j_2}(t_1,t_2)$ capturing the dependence between failure times of the two individuals with failure occurring due to a particular cause. For this purpose, we introduce frailty random variable(s) possibly depending on the cause of failure. Let us write, for $k=1,2$, the frailty vector $\boldsymbol{\epsilon^{(k)}} = (\epsilon^{(k)}_{1},\cdots,\epsilon^{(k)}_{L_k})$, 
where $\epsilon^{(k)}_{j}$ denotes the frailty variable, assumed to be positive, corresponding to  the $j$th cause of failure for the $k$th individual, for $j = 1,\cdots,L_k$. It is assumed that dependence between the two failure times $(T_1,T_2)$ with the corresponding causes $(J_1,J_2)$ is described solely through the random frailty terms $\boldsymbol{\epsilon^{(1)}}$ and $\boldsymbol{\epsilon^{(2)}}$. This means, conditional on ($\boldsymbol{\epsilon^{(1)}}, \boldsymbol{\epsilon^{(2)}})$, there is independence between $(T_1,J_1)$ and $(T_2,J_2)$. \\

Also, for $k=1,2$, let us introduce the cause-specific hazard 
$\lambda^{(k)}_{j}(t_k | \boldsymbol{\epsilon^{(k)}})$ for the failure of the $k$th individual at time $T_{k} = t_k$ due to cause $j$, conditional on the frailty vector $\boldsymbol{\epsilon^{(k)}}$, defined as 
$$\lambda^{(k)}_{j}(t_k | \boldsymbol{\epsilon^{(k)}}) = \lim\limits_{\delta \to 0+}\frac{P[t_{k} < T_{k} < t_{k} + \delta,J_{k} = j|T_{k} \geq t_{k},\boldsymbol{\epsilon^{(k)}}]}{\delta},$$
for $j = 1,\cdots,L_{k}$.   
As in the frailty literature (Wienke, 2010), we consider the multiplicative frailty model for $\lambda^{(k)}_{j}(t_k | \boldsymbol{\epsilon^{(k)}})$ given by 
\begin{equation}\label{EqGeneralmodel}
    \lambda^{(k)}_{j}(t_k | \boldsymbol{\epsilon^{(k)}}) = h^{(k)}_{0j}(t_k)\epsilon^{(k)}_{j}, 
\end{equation}
for $t_{k} > 0, j= 1,\cdots,L_{k}$ and 
$k=1,2$, where $h^{(k)}_{0j}(t_k)$ is the $j$th baseline cause-specific hazard at time $t_k$ for the $k$th individual. We assume, for each model, the existence of finite frailty mean. \\

The expression for the marginal survival function of the $k$th individual,  conditional on the frailty vector $\boldsymbol{\epsilon^{(k)}}$, is
\begin{align*}
    S^{(k)}(t_k | \boldsymbol{\epsilon^{(k)}}) 
    &= \exp{\Bigg[-\int\limits_{0}^{t_k}\sum\limits_{j= 1}^{L_k}\lambda^{(k)}_{j}(u_k | \boldsymbol{\epsilon^{(k)}})du_k\Bigg]}\\
    &= \exp{\Bigg[-\int\limits_{0}^{t_k}\sum\limits_{j = 1}^{L_k}h^{(k)}_{0j}(u_k)  \epsilon^{(k)}_{j}du_k\Bigg]} \\
    &= \exp{\Bigg[-\sum\limits_{j= 1}^{L_k}H^{(k)}_{0j}(t_k)  \epsilon^{(k)}_{j}\Bigg]},
\end{align*}
where $H^{(k)}_{0j}(t_k) = \int\limits_{0}^{t_{k}}h^{(k)}_{0j}(u)du$, for $k = 1,2$. Thus, the expression for $j$th sub-distribution function of  the $k$th individual, conditional on frailty vector $\boldsymbol{\epsilon^{(k)}}$, is 
\begin{align*}
   F_j^{(k)}(t_{k} | \boldsymbol{\epsilon^{(k)}})  
   &= \int\limits_{0}^{t_k}\lambda^{(k)}_{j}(u_k | \boldsymbol{\epsilon^{(k)}})S^{(k)}(u_k | \boldsymbol{\epsilon^{(k)}})du_k\\
   &= \int\limits_{0}^{t_k}h^{(k)}_{0j}(u_k)  \epsilon^{(k)}_{j}\exp{\Bigg[-\sum\limits_{j'= 1}^{L_k}H^{(k)}_{0j'}(u_k)  \epsilon^{(k)}_{j'}\Bigg]}du_k 
\end{align*}
for all $j = 1,\cdots,L_{k}$ and $k = 1,2$.\\

Therefore, using conditional independence given frailty, the joint sub-distribution function $F_{j_1j_2}(t_1,t_2)$ under the general model (\ref{EqGeneralmodel}) is
\begin{align}\label{eq:JointModel}
    F_{j_1j_2}(t_1,t_2) 
    &= \mathbb{E}\bigg[F^{(1)}_{j_1}(t_{1} | \boldsymbol{{\epsilon}^{(1)}})F^{(2)}_{j_2}(t_{2} | \boldsymbol{{\epsilon}^{(2)}})\bigg]\nonumber\\
    &= \mathbb{E}\Bigg[\int\limits_{0}^{t_1}\int\limits_{0}^{t_2} \prod_{k=1}^2 \left(h^{(k)}_{0j_{k}}(u_k) \epsilon^{(k)}_{j_k}\right) \exp{\Bigg(-\sum_{k=1}^2 
    \sum\limits_{j= 1}^{L_k}H^{(k)}_{0j}(u_k)  \epsilon^{(k)}_{j}\Bigg)}du_2 du_1\Bigg]\nonumber\\
    &= \int\limits_{0}^{\infty}\cdots\int\limits_{0}^{\infty}\Bigg[\int\limits_{0}^{t_1}\int\limits_{0}^{t_2}
    \prod_{k=1}^2 \left(h^{(k)}_{0j_{k}}(u_k) \epsilon^{(k)}_{j_k}\right) \exp{\Bigg(-\sum_{k=1}^2 
    \sum\limits_{j= 1}^{L_k}H^{(k)}_{0j}(u_k)  \epsilon^{(k)}_{j}\Bigg)}du_2 du_1\Bigg] \nonumber \\
    &\qquad\qquad\qquad\qquad\qquad\times g(\boldsymbol{\epsilon^{(1)}},\boldsymbol{\epsilon^{(2)}})d\boldsymbol{\epsilon^{(1)}} d\boldsymbol{\epsilon^{(2)}}, 
\end{align}
for all $j_{k} = 1,\cdots,L_{k}$ and $k = 1,2$, where  $g(\boldsymbol{\epsilon^{(1)}},\boldsymbol{\epsilon^{(2)}})$ is the joint density of the random frailty vector $(\boldsymbol{\epsilon^{(1)}},\boldsymbol{\epsilon^{(2)}})$. \\

Let us define the classes $\mathcal{H}_{k}$ of baseline cause-specific hazards for the $k$th individuals in a pair as \\
$$\displaystyle{\mathcal{H}_{k} = \bigg\{\Big(h^{(k)}_{01}(x),\cdots,h^{(k)}_{0L_k}(x)\Big): x > 0,\ \lim\limits_{x \to 0+}H^{(k)}_{0j}(x) = 0,\ j = 1,\cdots,L_{k}, \mbox{ and } \lim\limits_{x \to \infty}\sum_{j=1}^{L_k}H^{(k)}_{0j}(x) = \infty\bigg\}}, $$
and let $\boldsymbol{h^{(k)}_0}\in\mathcal{H}_k$ denote the set of these arbitrary baseline cause-specific hazard functions for the $k$th individual, for $k=1,2$. 
We also define the class of the frailty distributions 
$$\mathcal{G} =  \bigg\{g\Big(\boldsymbol{\epsilon^{(1)}},\boldsymbol{\epsilon^{(2)}}\Big)\ge 0:  \int\limits_{0}^{\infty}\cdots\int\limits_{0}^{\infty}g\Big(\boldsymbol{\epsilon^{(1)}},\boldsymbol{\epsilon^{(2)}}\Big)d\boldsymbol{\epsilon^{(1)}} d\boldsymbol{\epsilon^{(2)}} = 1\bigg\}.$$ 

\begin{definition}
 Model (\ref{EqGeneralmodel}) is identifiable within the class $\mathcal{H}_1\times\mathcal{H}_2\times\mathcal{G}$ if, for some $\boldsymbol{h^{(k)}_0}, \boldsymbol{\tilde{h}^{(k)}_0}\in\mathcal{H}_k$ for $k=1,2$, and 
$g\Big(\boldsymbol{\epsilon^{(1)}},\boldsymbol{\epsilon^{(2)}}\Big),\widetilde{g}\Big(\boldsymbol{\epsilon^{(1)}},\boldsymbol{\epsilon^{(2)}}\Big) \, \in \, \mathcal{G}$, 
the equality 
$$F_{j_1j_2}(t_1,t_2) = \widetilde{F}_{j_1j_2}(t_1,t_2) \, \, \text{for all} \, \, t_1,t_2 > 0 \, \, \text{and} \, \, j_k = 1,\cdots,L_{k},k = 1,2,$$
where $F_{j_1j_2}(t_1,t_2)$ and $\widetilde{F}_{j_1j_2}(t_1,t_2)$ are obtained from (\ref{eq:JointModel}) using the above-mentioned two sets of cause-specific hazards and frailty vector, respectively, implies that \\
$$h^{(k)}_{0j}(x) = \widetilde{h}^{(k)}_{0j}(x),$$ 
for all $x \ge 0$, \, $j = 1,\cdots,L_{k}, \, k = 1,2$ and 
$$g\Big(\boldsymbol{\epsilon^{(1)}},\boldsymbol{\epsilon^{(2)}}\Big)=\widetilde{g}\Big(\boldsymbol{\epsilon^{(1)}},\boldsymbol{\epsilon^{(2)}}\Big),$$
for all possible values of $\Big(\boldsymbol{\epsilon^{(1)}},\boldsymbol{\epsilon^{(2)}}\Big)$. 
\end{definition}
 
As described above, the classes $\mathcal{H}_1$,$\mathcal{H}_2$ and $\mathcal{G}$ are non-parametric in the sense that the baseline cause-specific hazards in $\mathcal{H}_1$,$\mathcal{H}_2$ and the frailty distributions in $\mathcal{G}$ do not belong to any specific parametric families. In the following, we prove that, in this generality, the model (\ref{EqGeneralmodel}) is not identifiable. For this, note that the expression for joint sub-density function is 
\begin{align}\label{jointsubdensity}
    f_{j_1j_2}(t_1,t_2) &= \int\limits_{0}^{\infty}\cdots\int\limits_{0}^{\infty} \prod_{k=1}^2 \left(     
    h^{(k)}_{0j_k}(t_k) \epsilon^{(k)}_{j_k}\right)
    \exp{\bigg[-\sum_{k=1}^2 \sum\limits_{j = 1}^{L_k}H^{(k)}_{0j}(t_k)\epsilon^{(k)}_{j} 
    \bigg]} \times \nonumber \\
    &\qquad\qquad g(\epsilon^{(1)}_{1},\cdots,\epsilon^{(1)}_{L_1},\epsilon^{(2)}_{1},\cdots,\epsilon^{(2)}_{L_2}) \prod_{k=1}^2\prod_{j=1}^{L_k}d\epsilon^{(k)}_{j},
\end{align}
where $h^{(k)}_{0j_k}(t_k)\in \mathcal{H}_{k}$, for  
 $t_k>0,\ j_{k} = 1,\cdots,L_{k}$ and $k = 1,2$. 
Let us write $$\tilde{h}^{(1)}_{0j_1}(t_1)= c_1h^{(1)}_{0j_1}(t_1) \quad\mbox{ and }\quad \tilde{h}^{(2)}_{0j_2}(t_2)=\frac{1}{c_2}h^{(2)}_{0j_2} (t_2),$$
for all $t_1,t_2>0$, $j_k = 1,\cdots,L_k$ and $k=1,2$, where $c_1$ and $c_2$ are some positive constants. It can be easily verified that the $\tilde{h}^{(k)}_{0j_k}(t_k)$, for $k=1,2$, are also valid baseline cause-specific hazard functions, thus belonging to $\mathcal{H}_1$ and $\mathcal{H}_2$, respectively. \\

Also, let us consider the frailty distribution given by $$\tilde{g}\Big(\boldsymbol{\epsilon^{(1)}},\boldsymbol{\epsilon^{(2)}}\Big)= \frac{c_1^{L_1}}{c_2^{L_2}} g\Big(c_1\boldsymbol{\epsilon^{(1)}},\frac{1}{c_2}\boldsymbol{\epsilon^{(2)}}\Big).$$
Since $g \in \mathcal{G}$, it is an easy exercise to check that $\tilde{g}$ also belongs to $\mathcal{G}$. Using the transformation of variables, 
$\boldsymbol{\delta^{(1)}}=c_1\boldsymbol{\epsilon^{(1)}}$ and $\boldsymbol{\delta^{(2)}}=\frac{1}{c_2}\boldsymbol{\epsilon^{(2)}},$ 
one can check that 
$\int_0^{\infty} \cdots \int_0^{\infty}  \tilde{g}\Big(\boldsymbol{\epsilon^{(1)}},\boldsymbol{\epsilon^{(2)}}\Big) d\boldsymbol{\epsilon^{(1)}} d\boldsymbol{\epsilon^{(2)}}=1.$ \\

With this choice of $\tilde{h}^{(k)}_{0j_k}(t_k)$, for $k=1,2$, and $\tilde{g}\Big(\boldsymbol{\epsilon^{(1)}},\boldsymbol{\epsilon^{(2)}}\Big)$, we have the corresponding joint sub-density function given by  
\begin{align*}
    \tilde{f}_{j_1j_2}(t_1,t_2) &= \int\limits_{0}^{\infty}\cdots\int\limits_{0}^{\infty}c_1h^{(1)}_{0j_1}(t_1)\frac{1}{c_2}h^{(2)}_{0j_2}(t_2) \epsilon^{(1)}_{j_1}\epsilon^{(2)}_{j_2}\exp{\bigg[-\sum\limits_{j = 1}^{L_1}c_1H^{(1)}_{0j}(t_1)\epsilon^{(1)}_{j}-\sum\limits_{j = 1}^{L_2}\frac{1}{c_2}H^{(2)}_{0j}(t_2)\epsilon^{(2)}_{j}\bigg]}\nonumber\\
    &\qquad\times \frac{c_1^{L_1}}{c_2^{L_2}} g(c_1\epsilon^{(1)}_{1},\cdots,c_1\epsilon^{(1)}_{L_1}, \frac{1}{c_2}\epsilon^{(2)}_{1},\cdots,\frac{1}{c_2}\epsilon^{(2)}_{L_2})\prod_{j = 1}^{L_1}d\epsilon^{(1)}_{j}\prod_{j = 1}^{L_2}d\epsilon^{(2)}_{j} \nonumber \\     &=\int\limits_{0}^{\infty}\cdots\int\limits_{0}^{\infty} \prod_{k=1}^2 \left(
    h^{(k)}_{0j_k}(t_k)\delta^{(k)}_{j_k}\right)
    \exp{\bigg[-\sum_{k=1}^2 \sum\limits_{j = 1}^{L_k}H^{(k)}_{0j}(t_k)\delta^{(k)}_{j} \bigg]}\nonumber\\
    &\qquad\times     
    g(\delta^{(1)}_{1},\cdots,\delta^{(1)}_{L_1}, \delta^{(2)}_{1},\cdots,\delta^{(2)}_{L_2}) \prod_{k=1}^2\prod_{j = 1}^{L_k}d\delta^{(k)}_{j},\nonumber\\ 
\end{align*}
which is same as $f_{j_1,j_2}(t_1,t_2)$ in (\ref{jointsubdensity}) for $t_k>0,\ j_{k} = 1,\cdots,L_{k}$ and $k = 1,2$. Therefore, when both the baseline cause-specific hazards and the frailty distribution are unspecified, the model (\ref{EqGeneralmodel}) is not identifiable. \\ 

In the following sections, we investigate the identifiability issues for the model (\ref{EqGeneralmodel}) with non-parametric baseline cause-specific hazards $h^{(k)}_{0j}(t_k)$'s, $j = 1,\cdots,L_{k}$, belonging to family $\mathcal{H}_k$, for $k=1,2$, but the frailty distribution given by the joint density $g(\boldsymbol{\epsilon^{(1)}},\boldsymbol{\epsilon^{(2)}})$ belonging to some parametric family with associated parameter vector $\boldsymbol{\theta}$. We define this parametric family by the set $\boldsymbol{\Theta}$ of all $\boldsymbol{\theta}$. In the next section, we give an example of a parametric family of frailty distributions that leads to non-identifiability of the model (\ref{EqGeneralmodel}) to show that parametric frailty does not always lead to identifiability. 
The subsequent sections consider different types of frailty distributions belonging to some Gamma family, since Gamma frailty seems to be commonly used by many authors for its convenience (Yashin et al., 1995; Iachine, 2004; Duchateau, 2008; Wienke, 2010 etc). 

\section{A non-identifiable model}

To work with a parametric family of frailty distribution, the joint sub-distribution function $F_{j_1j_2}(t_1,t_2)$ in (\ref{eq:JointModel}) is written as $F_{j_1j_2}(t_1,t_2;\boldsymbol{\theta})$ to indicate the dependence on the parameter vector $\boldsymbol{\theta}\in \boldsymbol{\Theta}$. 

\begin{definition}
    We say that the model $F_{j_1j_2}(t_1,t_2;\boldsymbol{\theta})$ is identifiable within the class $\mathcal{H}_{1}\times\mathcal{H}_{2}\times  \boldsymbol{\Theta}$ if, for some $\boldsymbol{h^{(k)}_0},\boldsymbol{\tilde{h}^{(k)}_0} \in \mathcal{H}_{k}$, for $k = 1,2$, and $\boldsymbol{\theta},\boldsymbol{\tilde{\theta}} \in \boldsymbol{\Theta}$, the equality 
      $$F_{j_1j_2}(t_1,t_2;\boldsymbol{\theta}) = F_{j_1j_2}(t_1,t_2;\boldsymbol{\tilde{\theta}}) \, \text{for all} \, \, t_{k} > 0, j_{k} = 1,\cdots,L_{k}, k=1,2,$$ 
      implies 
$$\boldsymbol{h^{(k)}_0}(x)=\boldsymbol{\tilde{h}^{(k)}_0}(x)\, \,\text{and} \, \,\boldsymbol{\theta} = \boldsymbol{\tilde{\theta}}\, \, \text{for all} \, \, x \ge 0, k = 1,2. $$
\end{definition} 

Let us consider the model 
\begin{equation}\label{Nonidentifiable}
  \lambda^{(k)}_{j}(t_{k} | \boldsymbol{\epsilon}^{(k)})  = h^{(k)}_{0}(t_k)  \epsilon^{(k)}_{j}, 
\end{equation}
for $t_{k} > 0,\ j = 1,\cdots,L_k$ and $k=1,2$, 
where $h^{(k)}_{0}(t_k)$ is the arbitrary $j$th baseline cause-specific hazard function for the $k$th individual, being the same for all $j = 1,\cdots,L_k$. Therefore, the class $\mathcal{H}_{k}$ now reduces to 
$$\displaystyle{\mathcal{H}_{k} = \bigg\{h^{(k)}_{0}(x): h^{(k)}_{0}(x) \ge 0, \int_{0}^{\infty}h^{(k)}_{0}(x)dx = \infty\bigg\}},$$
for $k=1,2$. 
Let us define the frailty variable $\epsilon^{(k)}_{j}$ in the general model (\ref{EqGeneralmodel}) as 
$$\epsilon^{(k)}_{j}=a\eta^{(k)}_{j},$$
for $j=1,\cdots,L_k,\ k=1,2$, where $a$ is a positive random variable, which acts as the shared frailty between the two individuals in a pair, and $\eta^{(k)}_{j}$ is $j$-th cause-specific frailty for the $k$-th individual. We assume that $\boldsymbol{\eta^{(k)}} = (\eta^{(k)}_{1},\cdots,\eta^{(k)}_{L_k}) \sim \text{Dir}(\boldsymbol{\alpha^{(k)}})$, for $k=1,2$, where 
Dir$(\boldsymbol{\alpha^{(k)}})$ denotes a Dirichlet distribution with parameter vector $\boldsymbol{\alpha^{(k)}}=(\alpha^{(k)}_{1},\cdots,\alpha^{(k)}_{L_k})$, with $\alpha^{(k)}_{j} > 0,\ j = 1,2,\cdots,L_{k}$, for $k=1,2$. 
While we assume independence between the two frailty vectors  $\boldsymbol{\eta^{(1)}}$ and $\boldsymbol{\eta^{(2)}}$, these are also assumed to be independent of the shared frailty $a$. Also, $a$ is assumed to follow a $\text{Gamma}(\frac{1}{\sigma^2},\frac{1}{\sigma^2})$ distribution, for some $\sigma > 0$, having density 
\begin{equation}\label{gammadensity}
g(\epsilon;\sigma^2)=\frac{1}{\sigma^{\frac{2}{\sigma^2}}\Gamma(\frac{1}{\sigma^2})} e^{-\frac{\epsilon}{\sigma^2}}
\epsilon^{\frac{1}{\sigma^2} - 1},
\end{equation}
so that $E[a]=1$. So, the parameter space $\boldsymbol{\Theta}$ is given by 
$$\boldsymbol{\Theta} = \big\{\boldsymbol{\theta}= (\boldsymbol{\alpha^{(1)}},\boldsymbol{\alpha^{(2)}},\sigma):\boldsymbol{\alpha^{(k)}} \in \mathbb{R}^{L_k}_{+}, \, \, \text{for}\, \, k = 1,2, \, \, \text{and} \, \, \sigma > 0\big\}.$$

The joint frailty density function $g(\boldsymbol{\epsilon^{(1)}},\boldsymbol{\epsilon^{(2)}},\boldsymbol{\theta})$ can be derived from the above description. The choice of Dirichlet distribution for cause-specific frailty vector $\boldsymbol{\eta^{(k)}}$ facilitates the derivation because we have  $\sum\limits_{j = 1}^{L_k}\eta^{(k)}_{j} = 1$ almost surely, for $k = 1,2$. Also, the Gamma distribution for the shared frailty variable $a$ makes the computation for joint sub-distribution function easier. 

Note that $\eta^{(k)}_{j}$ marginally follows the $\text{Beta}\Big(\alpha^{(k)}_{j},\sum\limits_{j'\ne j}\alpha^{(k)}_{j'}\Big)$ distribution. Therefore, we have 
$$\displaystyle{\mathbb{E}(\eta^{(k)}_{j}}) = \frac{\alpha^{(k)}_{j}}{\sum\limits_{j' = 1}^{L_k}\alpha^{(k)}_{j'}},$$ 
for $j = 1,\cdots,L_{k}$ and $k = 1,2$. Also, since $a$ is independent of the $\eta^{(k)}_{j}$'s and $\mathbb{E}(a)=1$, we have $\mathbb{E}(\epsilon^{(k)}_{j})=\mathbb{E}(\eta^{(k)}_{j})$. \\
The survival function for the $k$th individual, conditional on frailty variables $a,\boldsymbol{\eta^{(k)}}$, is given by 
\begin{align*}
    S^{(k)}(t | a,\boldsymbol{\eta^{(k)}}) 
    &= \exp\Bigg[-\int\limits_{0}^{t}\sum\limits_{j = 1}^{L_k}\lambda^{(k)}_{j}(u | a,\boldsymbol{\eta^{(k)}})du\Bigg]\\
    &= \exp\Bigg[-\int\limits_{0}^{t}a h^{(k)}_{0}(u)\sum\limits_{j = 1}^{L_k} \eta^{(k)}_{j}du \Bigg]\\
    &= \exp\Bigg[-\int\limits_{0}^{t}a h^{(k)}_{0}(u)du\Bigg]\\
&= \exp[-a H^{(k)}_{0}(t)], 
\end{align*}
where $H^{(k)}_{0}(t) = \int_{0}^{t}h^{(k)}_{0}(u)du$, for $k = 1,2$. Hence, the $j$th sub-distribution function for the $k$th individual, conditional on the frailty variables $a,\boldsymbol{\eta^{(k)}}$, is 
\begin{align*}
    F^{(k)}_{j}(t_k | a,\boldsymbol{\eta^{(k)}}) 
&= \int\limits_{0}^{t_k}\lambda^{(k)}_{j}(u_k | a,\boldsymbol{\eta^{(k)}})S^{(k)}(u_k | a,\boldsymbol{\eta^{(k)}})du_k\\
&= \int\limits_{0}^{t_k}h^{(k)}_{0}(u_k)a\eta^{(k)}_{j}\exp{\{-a H^{(k)}_{0}(u_k) \}}du_k\\
&= \eta^{(k)}_{j} \int\limits_{0}^{a H^{(k)}_{0}(t_k)} e^{-z} dz\\
&= \eta^{(k)}_{j}\big(1 - \exp{\{-a H^{(k)}_{0}(t_k)\}}\big), 
\end{align*}
for $t_k>0,\ j = 1,\cdots,L_{k}$ and $k = 1,2$. Therefore, the unconditional joint sub-distribution function $F_{j_{1}j_{2}}(t_1,t_2;\boldsymbol{\theta})$ is derived, using conditional independence given $a$,$\boldsymbol{\eta^{(1)}}$ and $\boldsymbol{\eta^{(2)}}$, as  
\begin{align*}
    F_{j_{1}j_{2}}(t_1,t_2;\boldsymbol{\theta}) &= \mathbb{E}\bigg[F^{(1)}_{j_{1}}(t_1 | a,\boldsymbol{\eta^{(1)}}) F^{(2)}_{j_{2}}(t_2 | a,\boldsymbol{\eta^{(2)}})\bigg] \\
    &= \mathbb{E}\bigg[\eta^{(1)}_{j_{1}}\eta^{(2)}_{j_{2}}\Big(1 - \exp{\big\{-a H^{(1)}_{0}(t_1)\big\}}\Big)\Big(1 - \exp{\big\{-a H^{(2)}_{0}(t_2)\big\}}\Big)\bigg]\\
    &= \left[\prod_{k=1}^2 \mathbb{E}(\eta^{(k)}_{j_{k}})\right] \times     
    \mathbb{E}\bigg[(1 - \exp{\big\{-a H^{(1)}_{0}(t_1)\big\}})(1 - \exp{\big\{-a H^{(2)}_{0}(t_2)\big\}})\bigg]\\
    &= \left[\prod_{k=1}^2 \mathbb{E}(\eta^{(k)}_{j_{k}})\right] \times
    \bigg[1 - \sum\limits_{k=1}^{2}\mathbb{E}\Big(\exp{\big\{-a H^{(k)}_{0}(t_k)\big\}}\Big) + \mathbb{E}\Big(\exp{\big\{-a\sum\limits_{k=1}^{2}H^{(k)}_{0}(t_k)\big\}}\Big)\bigg].
\end{align*}
Note that
  \begin{align*}
    \mathbb{E}\bigg[\exp{\Big[-aH^{(k)}_{0}(t_k)\Big]}\bigg] &= \int\limits_{0}^{\infty}\exp{\Big\{- aH^{(k)}_{0}(t_k)\Big\}}\frac{\exp{\left(-\frac{a}{\sigma^2}\right)} a^{\frac{1}{\sigma^2} - 1}}{\sigma^{\frac{2}{\sigma^2}}\Gamma(\frac{1}{\sigma^2})}da\\
    &= \frac{1}{\sigma^{\frac{2}{\sigma^2}}\Gamma(\frac{1}{\sigma^2})} \int\limits_{0}^{\infty}\exp{\Big[- a\Big(H^{(k)}_{0}(t_k) + \frac{1}{\sigma^2}\Big)\Big]}a^{\frac{1}{\sigma^2} - 1}da\\
    &= \frac{1}{\sigma^{\frac{2}{\sigma^2}}\Gamma(\frac{1}{\sigma^2})} \times \frac{\Gamma(\frac{1}{\sigma^2})}{\Big(H^{(k)}_{0}(t_k) + \frac{1}{\sigma^2}\Big)^{\frac{1}{\sigma^2}}}\\
    &= \frac{1}{\Big(1 + \sigma^2 H^{(k)}_{0}(t_k)\Big)^{\frac{1}{\sigma^2}}}, 
\end{align*}
for $t_{k} > 0$ and $k = 1,2$. Hence, we have 
$$F_{j_1j_2}(t_1,t_2;\boldsymbol{\theta}) = 
\left(\prod_{k=1}^2 \frac{\alpha^{(k)}_{j_{k}}}{\sum\limits_{j = 1}^{L_k}\alpha^{(k)}_{j}}\right)\times 
\bigg[1 - \sum\limits_{k=1}^{2}\Big(1 + \sigma^2 H^{(k)}_{0}(t_k)\Big)^{-\frac{1}{\sigma^2}} + \Big(1 + \sigma^2 \sum\limits_{k=1}^{2} H^{(k)}_{0}(t_k)\Big)^{-\frac{1}{\sigma^2}}\bigg],$$
for all $t_{k} > 0,\ j_{k} = 1,\cdots,L_{k}$ and $k = 1,2$.\\ 

Let us consider 
$$\boldsymbol{\tilde{\theta}} = (\tilde{\alpha}^{(1)}_{1},\cdots,\tilde{\alpha}^{(1)}_{L_1},\tilde{\alpha}^{(2)}_{1},\cdots,\tilde{\alpha}^{(2)}_{L_2},\tilde{\sigma}),$$ 
where $\tilde{\alpha}^{(k)}_{j} = c_{k} \alpha^{(k)}_{j}$, for some $c_k>0$ and $j = 1,\cdots,L_{k},\ k=1,2$; also, $\tilde{\sigma}=\sigma$. 
This implies $\sum\limits_{j = 1}^{L_k}\tilde{\alpha}^{(k)}_{j} = c_k\sum\limits_{j = 1}^{L_k}\alpha^{(k)}_{j}$, for $k=1,2$. We also take $\tilde{h}^{(k)}_{0}(x) = h^{(k)}_{0}(x)$, for all $x > 0$ and  $k = 1,2$. Then, the joint sub-distribution function $\tilde{F}_{j_1j_2}(t_1,t_2;\boldsymbol{\tilde{\theta}})$ simplifies to 
\begin{align*}
    \tilde{F}_{j_1j_2}(t_1,t_2;\boldsymbol{\tilde{\theta}})  &= 
    \left(\prod_{k=1}^2 \frac{\tilde{\alpha}^{(k)}_{j_{k}}}{\sum\limits_{j = 1}^{L_k}\tilde{\alpha}^{(k)}_{j}}\right)\times 
    \Bigg[1 - \sum\limits_{k=1}^{2}\bigg(1 + \sigma^2 H^{(k)}_{0}(t_k)\bigg)^{-\frac{1}{\sigma^2}} + \bigg(1 + \sigma^2 \sum\limits_{k=1}^{2} H^{(k)}_{0}(t_k)\bigg)^{-\frac{1}{\sigma^2}}\Bigg]\\
   &= \left(\prod_{k=1}^2 \frac{\alpha^{(k)}_{j_{k}}}{\sum\limits_{j = 1}^{L_k}\alpha^{(k)}_{j}}\right)\times    
 \Bigg[1 - \sum\limits_{k=1}^{2}\bigg(1 + \sigma^2 H^{(k)}_{0}(t_k)\bigg)^{-\frac{1}{\sigma^2}} + \bigg(1 + \sigma^2 \sum\limits_{k=1}^{2} H^{(k)}_{0}(t_k)\bigg)^{-\frac{1}{\sigma^2}}\Bigg]\\
   &= F_{j_1j_2}(t_1,t_2;\boldsymbol{\theta}), 
\end{align*}
for all $t_{k} > 0 $ and $j_{k} = 1,\cdots,L_{k}$ and $k=1,2$. Hence, this model (\ref{Nonidentifiable}) is non-identifiable.\\

Nevertheless, in the following four sections, we consider four different types of Gamma frailty variables: (1) shared Gamma frailty, (2) correlated Gamma frailty, (3) shared cause-specific Gamma frailty and (4) correlated cause-specific Gamma frailty, respectively.  

\section{Shared Gamma frailty model}

Shared Gamma frailty is the simplest model involving one common frailty variable which is shared between the two individuals in a pair regardless of failure type. It captures the common pair or cluster effect (\cite{wienke2010frailty}), thereby inducing association between the two individuals. An interesting feature of this model is that one can consider two different sets of competing risks corresponding to the two individuals in a pair. For example, while studying multi-organ failure, we may encounter different sets of failure types corresponding to different organs. In such situations, shared frailty captures possible dependence due to a common body effect. 
In this case, the model (\ref{EqGeneralmodel}) reduces to 
\begin{align}\label{NonparametrichazardsharedGamma}
   \lambda^{(k)}_{j}(t_k | \epsilon) &= h^{(k)}_{0j}(t_k)\epsilon 
\end{align}
where $\epsilon_j^{(k)}=\epsilon$ is the common shared frailty variable, for $j=1,\cdots,L_k$ and $k=1,2$. The shared frailty $\epsilon$ is assumed to follow a $\text{Gamma}(\frac{1}{\sigma^2},\frac{1}{\sigma^2})$ distribution with density given by (\ref{gammadensity}), for some $\sigma > 0$. So, the parameter space is 
$$\boldsymbol{\Theta} = \big\{\sigma:\sigma > 0\big\}.$$ 

Following the same technique as before, the marginal survival function of the $k$th individual, conditional on the shared frailty $\epsilon$, is 
    $$S^{(k)}(t;\mathbf{h^{(k)}_0} | \epsilon) = 
    \exp{\bigg[-\epsilon \sum_{j=1}^{L_k}H^{(k)}_{0j}(t)\bigg]}=
    \exp{\bigg[-\epsilon H^{(k)}_{0}(t)\bigg]}, $$
where $H^{(k)}_{0}(t)=\sum_{j=1}^{L_k}H^{(k)}_{0j}(t)$, unlike in Section 3, 
for $k=1,2$. Also, the marginal sub-distribution function due to cause $j$ for the $k$th individual, conditional on the shared frailty  $\epsilon$, is 
\begin{align*}
    F^{(k)}_{j}(t_k;\mathbf{h^{(k)}_0} | {\epsilon}) 
    &= \int\limits_{0}^{t_k}h^{(k)}_{0j}(u_k)\epsilon\exp{\big\{-\epsilon H^{(k)}_{0}(u_k) \big\}}du_k,
\end{align*}
for $j = 1,2,\cdots,L_k$ and $k = 1,2$. Therefore, as before, the unconditional joint sub-distribution function $F_{j_1j_2}(t_1,t_2;\mathbf{h^{(1)}_0},\mathbf{h^{(2)}_0},\sigma^2)$ is written as 
\begin{align*}
    F_{j_1j_2}(t_1,t_2;\mathbf{h^{(1)}_0},\mathbf{h^{(2)}_0},\sigma^2) &= \mathbb{E}\bigg[F^{(1)}_{j_1}(t_1;\mathbf{h^{(1)}_0} | \epsilon)F^{(2)}_{j_2}(t_2;\mathbf{h^{(2)}_0}| \epsilon)\bigg]\\
    &= \int\limits_{0}^{\infty}\bigg[\int\limits_{0}^{t_1}\int\limits_{0}^{t_2} \left(\prod_{k=1}^2 h^{(k)}_{0j_k}(u_k)\right)\epsilon^{2} \times\exp{\big(-\epsilon \sum_{k=1}^2 
    H^{(k)}_{0}(u_k)\big)}du_2du_1\bigg]g(\epsilon;\sigma^2)d\epsilon, 
\end{align*}
where $g(\epsilon;\sigma^2)$ is the density of the $\text{Gamma}(\frac{1}{\sigma^2},\frac{1}{\sigma^2})$ distribution given by 
(\ref{gammadensity}). 
This simplifies the unconditional joint sub-distribution function $F_{j_1j_2}(t_1,t_2;\mathbf{h^{(1)}_0},\mathbf{h^{(2)}_0},\sigma^2)$ 
as 
\begin{align*}
&=\int\limits_{0}^{\infty}\Bigg[\int\limits_{0}^{t_1}\int\limits_{0}^{t_2} \left(\prod_{k=1}^2 h^{(k)}_{0j_k}(u_k)\right)\epsilon^{2} \times\exp{\big(-\epsilon \sum_{k=1}^2 
    H^{(k)}_{0}(u_k)\big)}du_2du_1\bigg] \frac{\exp{\left(-\frac{\epsilon}{\sigma^2}\right)} \epsilon^{\frac{1}{\sigma^2} - 1}}{\sigma^{\frac{2}{\sigma^2}}\Gamma(\frac{1}{\sigma^2})}d\epsilon  \\ 
&=\int\limits_{0}^{t_1}\int\limits_{0}^{t_2}\Bigg[\int\limits_{0}^{\infty} \frac{1}{\sigma^{\frac{2}{\sigma^2}}\Gamma(\frac{1}{\sigma^2})}\exp{\Big(-\epsilon \big(\sum_{k=1}^2 
    H^{(k)}_{0}(u_k)+\frac{1}{\sigma^2}\big)\Big)}
\epsilon^{2 + \frac{1}{\sigma^2} - 1} d\epsilon\Bigg]\left(\prod_{k=1}^2 h^{(k)}_{0j_k}(u_k)\right) du_2du_1 \\
    &= \frac{1}{\sigma^{\frac{2}{\sigma^2}}\Gamma(\frac{1}{\sigma^2})}\int\limits_{0}^{t_1}\int\limits_{0}^{t_2} \left(\prod_{k=1}^2 h^{(k)}_{0j_k}(u_k)\right) 
    \frac{\Gamma(2 + \frac{1}{\sigma^2})}{\big(\sum_{k=1}^2 
    H^{(k)}_{0}(u_k)+\frac{1}{\sigma^2}\big)^{2 + \frac{1}{\sigma^2}}}du_2du_1 \\  
    &=\frac{1}{\sigma^{\frac{2}{\sigma^2}}\Gamma(\frac{1}{\sigma^2})}\int\limits_{0}^{t_1} \int\limits_{0}^{t_2} \frac{(1 + \frac{1}{\sigma^2}) \times \frac{1}{\sigma^2} \times \Gamma(\frac{1}{\sigma^2})\left(\prod_{k=1}^2 h^{(k)}_{0j_k}(u_k)\right)}{\big(\sum_{k=1}^2 
    H^{(k)}_{0}(u_k)+\frac{1}{\sigma^2}\big)^{2 + \frac{1}{\sigma^2}}}du_2du_1 \\   
    &=\int\limits_{0}^{t_1} \int\limits_{0}^{t_2} \frac{(1 + \sigma^2)\left(\prod\limits_{k=1}^2 h^{(k)}_{0j_k}(u_k)\right)}{\big(1+\sigma^2\sum\limits_{k=1}^2 
    H^{(k)}_{0}(u_k)\big)^{2 + \frac{1}{\sigma^2}}}du_2du_1, 
\end{align*}
for $t_k > 0,\ j_k = 1,\cdots,L_{k}$ and $k = 1,2$. Similarly, the unconditional sub-distribution function $F_{j}^{(k)}(t_k;\mathbf{h^{(k)}_0},\sigma^2)$ is obtained as 
\begin{align}\label{SubdistriSharedGamma}
&=\int\limits_{0}^{\infty}\Bigg[\int\limits_{0}^{t_k} h^{(k)}_{0j}(u_k)\epsilon 
\exp{\left\{-\epsilon H^{(k)}_{0}(u_k)\right\}}du_{k}\Bigg] \frac{\exp{\left(-\frac{\epsilon}{\sigma^2}\right) \epsilon^{\frac{1}{\sigma^2} - 1}}}{\sigma^{\frac{2}{\sigma^2}}\Gamma(\frac{1}{\sigma^2})}d\epsilon \nonumber \\ 
&= \int\limits_{0}^{t_k} \frac{h^{(k)}_{0j}(u_k)} {\big(1+\sigma^2 
    H^{(k)}_{0}(u_k)\big)^{1 + \frac{1}{\sigma^2}}}du_k, 
\end{align}
for $j=1,\cdots,L_k$ and $k=1,2$. The general definition of identifiability in Section 2 now reduces to the following. 
\begin{definition}
    The shared Gamma frailty model (\ref{NonparametrichazardsharedGamma}) is identifiable within the family $\mathcal{H}_{1}\times\mathcal{H}_{2}\times\Theta$ if, for some $\boldsymbol{h^{(k)}_0},\boldsymbol{\tilde{h}^{(k)}_0} \in \mathcal{H}_k$, for $k=1,2,$ and $\sigma,\ \tilde{\sigma} \in \boldsymbol{\Theta}$, the equality
$$F_{j_1j_2}(t_1,t_2;\boldsymbol{h^{(1)}_0},\boldsymbol{h^{(2)}_0},\sigma^2) = F_{j_1j_2}(t_1,t_2;\boldsymbol{\tilde{h}^{(1)}_0},\boldsymbol{\tilde{h}^{(2)}_0},\tilde{\sigma^2}),$$
for all $t_{k} > 0,\ j_{k} = 1,\cdots,L_{k}$ and $k = 1,2$, 
implies 
\begin{center}
$h^{(k)}_{0j}(x) = \tilde{h}^{(k)}_{0j}(x)$ for all $x > 0,\ j= 1,\cdots,L_{k},\ k = 1,2$ and $\sigma = \tilde{\sigma}.$  
\end{center}
\end{definition} 

\begin{theorem} 

The shared Gamma frailty model (\ref{NonparametrichazardsharedGamma}) with non-parametric baseline 
cause-specific hazard functions is identifiable  within the family  $\mathcal{H}_1\times\mathcal{H}_2\times\boldsymbol{\Theta}$, where $\boldsymbol{\Theta}= \big\{\sigma:\sigma > 0\big\}$, provided $H^{(k)}_{0j}(\cdot)$ is continuous and its inverse, denoted by $H^{-(k)}_{0j}(\cdot)$, exists for all $j=1,\cdots,L_k$ and $k=1,2$. 
\end{theorem}
\begin{proof}
From the equality of the joint of the sub-distribution functions as in {\bf Definition 2}, for all $t_k>0,\ j=1,\cdots,L_k$ and $k=1,2$, we get equality of the marginal sub-distribution functions $F_{j}^{(k)}(t_k;\mathbf{h^{(k)}_0},\sigma^2)$ and $F_{j}^{(k)}(t_k;\mathbf{\tilde{h}^{(k)}_0},\tilde{\sigma}^2)$, for all $t_k>0,\ j=1,\cdots,L_k$ and $k=1,2$. In particular, we have equality of the corresponding marginal sub-density functions obtained from  (\ref{SubdistriSharedGamma}), as given by 
\begin{equation*}
  \frac{h^{(k)}_{0j}(t_k)}{[1 + \sigma^{2} \sum_{j'=1}^{L_k} H^{(k)}_{0j'}(t_k)]^{1 + \frac{1}{\sigma^2}}} =  \frac{\tilde{h}^{(k)}_{0j}(t_k)}{[1 + \tilde{\sigma}^{2}\sum_{j'=1}^{L_k} \tilde{H}^{(k)}_{0j'}(t_k)]^{1 + \frac{1}{\tilde{\sigma}^2}}},
\end{equation*}
for all $t_k>0,\ j=1,\cdots,L_k$ and $k=1,2$. This equality can be written as 
\begin{equation*}
  \frac{\frac{d}{dt_k} H^{(k)}_{0j}(t_k)}{\frac{d}{dt_k} \tilde{H}^{(k)}_{0j}(t_k)}   
  =  \frac{[1 + \sigma^{2} \sum_{j'=1}^{L_k} H^{(k)}_{0j'}(t_k)]^{1 + \frac{1}{\sigma^2}}}{[1 + \tilde{\sigma}^{2}\sum_{j'=1}^{L_k} \tilde{H}^{(k)}_{0j'}(t_k)]^{1 + \frac{1}{\tilde{\sigma}^{2}}}}. 
\end{equation*} 
Now, writing $y=\sigma^{2} H^{(k)}_{0j}(t_k)$ and $z=\tilde{\sigma}^{2}\tilde{H}^{(k)}_{0j}(t_k)$, this equality reduces to 
\begin{equation*}
   \frac{\tilde{\sigma}^{2}}{\sigma^2}\times \frac{dy}{dz}=
  \frac{[1 + \sigma^{2} \sum_{j'(\ne j)=1}^{L_k} H^{(k)}_{0j'}(H^{-(k)}_{0j}(\frac{y}{\sigma^2}))+ y]^{1 + \frac{1}{\sigma^2}}}{[1 + \tilde{\sigma}^{2}\sum_{j'(\ne j)=1}^{L_k} \tilde{H}^{(k)}_{0j'}(\tilde{H}^{-(k)}_{0j}(\frac{z}{\tilde{\sigma}^{2}}))+z]^{1 + \frac{1}{\tilde{\sigma}^{2}}}}. 
\end{equation*}
This is equivalent to, after integration from 0 to some $w>0$, 
\begin{equation*}
\frac{\int\limits_{0}^{w} [1 + \sigma^{2} \sum_{j'(\ne j)=1}^{L_k} H^{(k)}_{0j'}(H^{-(k)}_{0j}(\frac{y}{\sigma^{2}}))+ y]^{-1 - \frac{1}{\sigma^2}} dy}{\int\limits_{0}^{w} [1 + \tilde{\sigma}^{2} \sum_{j'(\ne j)=1}^{L_k} \tilde{H}^{(k)}_{0j'}(\tilde{H}^{-(k)}_{0j}(\frac{z}{\tilde{\sigma}^{2}}))+ z]^{-1 - \frac{1}{\tilde{\sigma}^{2}}}dz} = \frac{\sigma^2}{\tilde{\sigma}^2}.
\end{equation*} 
Now, taking the limit as $w \to 0+$  and using the L'opital rule, we have 
\begin{equation*}
    \frac{\lim\limits_{w \to 0+}\Bigg[1 + \sigma^{2} \sum_{j'(\ne j)=1}^{L_k} H^{(k)}_{0j'}(H^{-(k)}_{0j}(\frac{w}{\sigma^{2}}))+ w\Bigg]^{-1 - \frac{1}{\sigma^2}}}{\lim\limits_{w \to 0+}\Bigg[1 + \tilde{\sigma}^{2} \sum_{j'(\ne j)=1}^{L_k} \tilde{H}^{(k)}_{0j'}(\tilde{H}^{-(1)}_{0j}(\frac{w}{\tilde{\sigma}^{2}}))+ w\Bigg]^{-1 - \frac{1}{\tilde{\sigma}^{2}}}}= \frac{\sigma^2}{\tilde{\sigma}^2}. 
\end{equation*}
Since the composite functions appearing in the above equality are also continuous functions, we have 
$$\lim\limits_{w \to 0+}
H^{(k)}_{0j'}(H^{-(k)}_{0j}(\frac{w}{\sigma^2}))= \lim\limits_{w \to 0+}\tilde{H}^{(1)}_{0j'}(\tilde{H}^{-(k)}_{0j}(\frac{w}{\tilde{\sigma}^2}))=0,$$
for $j\ne j'$. Hence, we get $\frac{\sigma^2}{\tilde{\sigma}^2}=1$, or $\sigma=\tilde{\sigma}$. \\
Now, equality of the joint sub-distribution functions also imply the equality of the marginal survival functions 
$S^{(k)}(t_k;\mathbf{h^{(k)}_0},\sigma^2)$ and $S^{(k)}(t_k;\boldsymbol{\tilde{h}^{(k)}_0},\tilde{\sigma}^2)$,  
for all $t_k>0$ and $k=1,2$. This gives, with $\sigma=\tilde{\sigma}$, 
$$\big[1 + \sigma^{2} H^{(k)}_{0}(t_k)\big]^{-\frac{1}{\sigma^2}} = \big[1 + \sigma^{2} \tilde{H}^{(k)}_{0}(t_k)\big]^{-\frac{1}{\sigma^{2}}},$$
for all $t_k > 0$, which implies 
    $$H^{(k)}_{0}(t_k) = \tilde{H}^{(k)}_{0}(t_k),$$
for all $t_k > 0$ and $k=1,2$. 

Again, from the equality of the sub-density functions  
$f^{(k)}_{j}(t_k;\mathbf{h^{(k)}_0},\sigma^2)$ and $f^{(k)}_{j}(t_k;\boldsymbol{\tilde{h}^{(k)}_0},\tilde{\sigma}^2)$ and using 
$\sigma=\tilde{\sigma}$ and $H^{(k)}_{0}(t_k) = \tilde{H}^{(k)}_{0}(t_k)$, for all $t_k > 0$, we get 
$$\frac{h^{(k)}_{0j}(t_k)}{\big[1 + \sigma^{2} H^{(k)}_{0}(t_k)\big]^{1 + \frac{1}{\sigma^2}}}  = \frac{\tilde{h}^{(k)}_{0j}(t_k)}{\big[1 + \sigma^{2} H^{(k)}_{0}(t_k)\big]^{1 + \frac{1}{\sigma^{2}}}}.$$
This implies 
$$h^{(k)}_{0j}(t_k) = \tilde{h}^{(k)}_{0j}(t_k),$$
for all $t_k > 0,\ j = 1,\cdots,L_k$ and $k=1,2$. \\ 
Hence, the shared Gamma frailty model (\ref{NonparametrichazardsharedGamma}) is identifiable. 
\end{proof}

\section{Correlated Gamma frailty model}

The correlated frailty model considers different, but correlated, frailty variables for the two individuals in a pair. These two frailty variables, denoted by $\epsilon^{(k)},\ k=1,2$, are only individual-specific, as in the shared frailty model, 
regardless of the cause of death. As a result, this model also 
can consider two different sets of competing risks corresponding to the two individuals in a pair. Considering the example of multi-organ failure of the previous section, correlated frailty model captures possible dependence due to a common body effect through the correlation between $\epsilon^{(1)}$ and $\epsilon^{(2)}$, while allowing for different frailty for different organs.

Formally, the correlated frailty model is given by 
\begin{equation}\label{NonparametrichazardcorrelatedGamma}
   \lambda^{(k)}_{j}(t_k | \epsilon^{(k)}) = h^{(k)}_{0j}(t_k) \epsilon^{(k)},
\end{equation}
for $t_{k} > 0,\  j = 1,\cdots,L_k$ and $k = 1,2$, where  $\epsilon^{(k)}$ is the frailty variable for the  $k$th individual. We assume that the frailty vector $(\epsilon^{(1)},\epsilon^{(2)})$ has a correlated Gamma distribution (See Yashin et al., 1995; Wienke, 2010), as described below. Write 
$\epsilon^{(k)} = \frac{\mu_0}{\mu_k}Y_0 + Y_k$, for $k=1,2$, where $Y_0 \sim \text{Gamma}(\kappa_0,\mu_0),Y_1 \sim \text{Gamma}(\kappa_1,\mu_1)$ and $Y_2 \sim \text{Gamma}(\kappa_2,\mu_2)$, with $\mu_0,\mu_1,\mu_2,\kappa_0,\kappa_1,\kappa_2 > 0$ satisfying $\mu_k = \kappa_0 + \kappa_k$, for $k=1,2$, and the random variables $Y_0,Y_1,Y_2$ are independent. This last condition is imposed to have expectation of $\epsilon^{(k)}$ equal to unity, as desired. 
Following Wienke (2010), $\epsilon^{(k)} \sim \text{Gamma}(\kappa_{0} + \kappa_{k},\mu_k)$, for $k = 1,2$. 
Therefore, $\mathbb{E}(\epsilon^{(k)}) = 1, \text{Var}(\epsilon^{(k)}) = 1/\mu_k$ and the covariance between $\epsilon^{(1)}$ and $\epsilon^{(2)}$ is given by 
$\text{Cov}(\epsilon^{(1)},\epsilon^{(2)}) = \text{Cov}(\frac{\mu_0}{\mu_1}Y_{0},\frac{\mu_0}{\mu_2}Y_{0}) = \frac{\mu^{2}_0}{\mu_{1}\mu_{2}}\times\frac{\kappa_0}{\mu^{2}_0} = \frac{\kappa_0}{(\kappa_{0} + \kappa_{1})(\kappa_{0} + \kappa_{2})}$.
 Clearly, there are three independent parameters in this correlated Gamma frailty model. We consider a re-parametrization of these parameters as 
$$\sigma_1=\frac{1}{\sqrt{\kappa_{0} + \kappa_{1}}},\ \sigma_2=\frac{1}{\sqrt{\kappa_{0} + \kappa_{2}}},\ \rho=\frac{\kappa_0}{\sqrt{\kappa_{0} + \kappa_{1}}\sqrt{\kappa_{0} + \kappa_{2}}},$$
so that $\sigma_{k}>0$ is the standard deviation of the frailty variable $\epsilon^{(k)}$, for $k=1,2$, and $\rho$ is the  correlation coefficient between $\epsilon^{(1)}$ and $\epsilon^{(2)}$, satisfying $0 < \rho < \min{(\frac{\sigma_1}{\sigma_2},\frac{\sigma_2}{\sigma_1})}$. Therefore, the frailty parameter space  $\boldsymbol{\Theta}$ is given by 
$$\boldsymbol{\Theta}=\big\{\boldsymbol{\theta}=(\sigma_1,\sigma_2,\rho): \sigma_{1} > 0, \sigma_{2} > 0,\  0 < \rho < \min{(\frac{\sigma_1}{\sigma_2},\frac{\sigma_2}{\sigma_1})}\big\}.$$\label{Theta_3}

As before, the survival function of the $k$th individual, conditional on  the frailty variable $\epsilon^{(k)}$, is given by 
$$S^{(k)}(t_k;\boldsymbol{h^{(k)}_{0}}| \epsilon^{(k)}) 
    = \exp{\Big[-\epsilon^{(k)}H^{(k)}_{0}(t_k)\Big]}= \exp{\Big[-
   H^{(k)}_{0}(t_k) (\frac{\mu_0}{\mu_k}Y_0 + Y_k)\Big]},$$
for $k = 1,2$. Also, the $j$th sub-distribution function of the $k$th individual, conditional on $\epsilon^{(k)}$, is 
$$ F^{(k)}_{j}(t_k;\boldsymbol{h^{(k)}_{0}} | \epsilon^{(k)}) = 
\int\limits_{0}^{t_k}h^{(k)}_{0j}(u_k)(\frac{\mu_0}{\mu_k}Y_0 + Y_k) \exp{\Big[-H^{(k)}_{0}(u_k)(\frac{\mu_0}{\mu_k}Y_0 + Y_k)\Big]}du,$$
 for $j = 1,\cdots,L_k$ and $k=1,2$. Then, $F_{j_1j_2}(t_1,t_2;\boldsymbol{h^{(1)}_{0}},\boldsymbol{h^{(2)}_{0}},\boldsymbol{\theta})$, the unconditional joint sub-distribution function, is given by 
$$\int\limits_{0}^{\infty}\int\limits_{0}^{\infty}\Bigg[\int\limits_{0}^{t_1}\int\limits_{0}^{t_2}\left(\prod_{k=1}^2 
h^{(k)}_{0j_k}(u_k)\epsilon^{(k)}\right)\times 
\exp{\Big(-\sum_{k=1}^2 
\epsilon^{(k)} H^{(k)}_{0}(u_k)\Big)}du_2 du_1\Bigg] g(\epsilon^{(1)},\epsilon^{(2)};\boldsymbol{\theta})
    d\epsilon^{(1)} d\epsilon^{(2)}\label{(***)},$$
for all $t_{k} > 0,\ j_{k} = 1,\cdots,L_{k}$ and $k=1,2$, where $g(\epsilon^{(1)},\epsilon^{(2)};\boldsymbol{\theta})$ denotes the joint density of $(\epsilon^{(1)},\epsilon^{(2)})$ following a correlated Gamma distribution, as described above. \\

From the expression of the $j$th sub-distribution function of the $k$th individual, conditional on $\epsilon^{(k)}$, or equivalently $Y_0$ and $Y_k$, we get the corresponding unconditional $j$th sub-distribution function, for $k=1,2$, as 
\begin{align*}
    F^{(k)}_{j}(t_k;\boldsymbol{h^{(k)}_{0}},\boldsymbol{\theta}) &= \int\limits_{0}^{\infty}\int\limits_{0}^{\infty}\Bigg[\int\limits_{0}^{t_k} h^{(k)}_{0j}(u_k)\bigg(\frac{\mu_0}{\mu_k}y_0 + y_k\bigg) \exp{\bigg( - H^{(k)}_{0}(u_k) \left(\frac{\mu_0}{\mu_k}y_0 + y_k\right)\bigg)}du_k\Bigg]\\
    &\qquad\times \prod_{l=0,k}\left( 
    \frac{\mu^{\kappa_l}_l}{\Gamma(\kappa_l)}\exp{[-\mu_{l}y_{l}]} y^{\kappa_{l} - 1}_{l}dy_{l}\right) \\
    &= \int\limits_{0}^{t_k}h^{(k)}_{0j}(u_k)\big[I_{1}(u_k) + I_{2}(u_k)\big]du_k \, \, \text{(say)}, 
    \end{align*}
where 
\begin{align*}
I_{1}(u_k) &= \frac{\mu^{\kappa_0}_0}{\Gamma(\kappa_0)}\frac{\mu^{\kappa_k}_k}{\Gamma(\kappa_k)}\frac{\mu_{0}}{\mu_{k}}\times\int\limits_{0}^{\infty}\int\limits_{0}^{\infty}\exp{\bigg[ - \bigg(\frac{H^{(k)}_{0}(u_k)}{\mu_k} + 1\bigg)\sum_{l=0,k}
\mu_{l}y_{l}\bigg]}y^{\kappa_{0}}_{0}y^{\kappa_{k} - 1}_{k}dy_{0}dy_{k}\\
&= \frac{\kappa_{0}}{\mu_k}\bigg(\frac{H^{(k)}_{0}(u_k)}{\mu_k} + 1\bigg)^{-\kappa_{0} - 1 - \kappa_{k}} 
=\frac{\kappa_{0}}{\mu_k}\bigg(\frac{H^{(k)}_{0}(u_k)}{\mu_k} + 1\bigg)^{- \mu_{k} - 1}
\end{align*}
and 
\begin{align*}
I_{2}(u_k) &= \left(\prod_{l=0,k}
\frac{\mu^{\kappa_l}_l}{\Gamma(\kappa_l)}\right)
\int\limits_{0}^{\infty}\int\limits_{0}^{\infty}\exp{\bigg[ - \bigg(\frac{H^{(k)}_{0}(u_k)}{\mu_k} + 1\bigg) \sum_{l=0,k}
\mu_{l}y_{l}\bigg]}y^{\kappa_{0} - 1}_{0}y^{\kappa_{k}}_{k}dy_{0}dy_{k}\\
&= \frac{\kappa_{k}}{\mu_k}\bigg(\frac{H^{(k)}_{0}(u_k)}{\mu_k} + 1\bigg)^{-\kappa_{0} - 1 - \kappa_{k}} 
= \frac{\kappa_{k}}{\mu_k}\bigg(\frac{H^{(k)}_{0}(u_k)}{\mu_k} + 1\bigg)^{- \mu_{k} - 1},
\end{align*}
for all $u_k > 0$. Hence,
\begin{align*}    
F^{(k)}_{j}(t_k;\boldsymbol{h^{(k)}_{0}},\boldsymbol{\theta}) &= 
\int\limits_{0}^{t_k}h^{(k)}_{0j}(u_k)\Bigg[\frac{\kappa_{0} + \kappa_{k}}{\mu_k}\bigg(\frac{H^{(k)}_{0}(u_k)}{\mu_k} + 1\bigg)^{- \mu_{k} - 1}\Bigg]du_k\\
&=\int\limits_{0}^{t_k}\frac{h^{(k)}_{0j}(u_k)du_k}{(1 + \sigma^{2}_{k}H^{(k)}_{0}(u_k))^{1 + \frac{1}{\sigma^{2}_k}}},
\end{align*}
for all $t_{k} > 0,\ j = 1,\cdots,L_{k}$ and $k=1,2$. Note that this unconditional $j$th sub-distribution function $F^{(k)}_{j}(t_k;\boldsymbol{h^{(k)}_{0}},\boldsymbol{\theta})$ has the same form as that in (\ref{SubdistriSharedGamma}) 
for the shared Gamma frailty model of the previous section with $\sigma^{2}_{k}$ replacing $\sigma^2$. Since this unconditional $j$th sub-distribution function depends on $\boldsymbol{\theta}$ only through $\sigma_k^2$, we write this as $F^{(k)}_{j}(t_k;\boldsymbol{h^{(k)}_{0}},\sigma_k^2)$. \\ 

\begin{definition}
The correlated Gamma frailty model (\ref{NonparametrichazardcorrelatedGamma}) is identifiable within the family $\mathcal{H}_{1}\times\mathcal{H}_2\times \boldsymbol{\Theta}$ if, for some $\boldsymbol{h^{(k)}_{0}},\boldsymbol{\tilde{h}^{(k)}_{0}} \in \mathcal{H}_{k}$, for $k=1,2$, and $\boldsymbol{\theta}, \boldsymbol{\tilde{\theta}}\in \boldsymbol{\Theta}$, where $\boldsymbol{\theta} = (\sigma_1,\sigma_2,\rho), \boldsymbol{\tilde{\theta}} = (\tilde{\sigma}_1,\tilde{\sigma}_2,\tilde{\rho})$, the equality of the joint sub-distribution functions 
$$F_{j_1j_2}(t_1,t_2;\boldsymbol{h^{(1)}_{0}},\boldsymbol{h^{(2)}_{0}},\boldsymbol{\theta}) = F_{j_1j_2}(t_1,t_2;\boldsymbol{\tilde{h}^{(1)}_{0}},\boldsymbol{\tilde{h}^{(2)}_{0}},\boldsymbol{\tilde{\theta}}),$$
for all $t_k>0,\ j_{k} = 1,\cdots,L_{k}$ and $k = 1,2$, implies 
$$h^{(k)}_{0j}(x) = \tilde{h}^{(k)}_{0j}(x) \, \, \text{for all} \, \, x>0,\ j= 1,\cdots,L_{k},\, \, \sigma_{k} = \tilde{\sigma}_{k}, \, \, \text{for} \, \, k = 1,2, \, \, \text{and} \, \,\rho = \tilde{\rho}.$$
\end{definition}

\begin{theorem}\label{Theorem 4.2}
    The correlated Gamma frailty model (\ref{NonparametrichazardcorrelatedGamma}) with non-parametric baseline cause-specific hazard functions is identifiable within the family  $\mathcal{H}_{1}\times\mathcal{H}_2\times\boldsymbol{\Theta}$, where $\boldsymbol{\Theta}=\big\{\boldsymbol{\theta}=(\sigma_1,\sigma_2,\rho): \sigma_{1} > 0, \sigma_{2} > 0,\  0 < \rho \leq \min{(\frac{\sigma_1}{\sigma_2},\frac{\sigma_2}{\sigma_1})}\big\}$, 
    provided $H^{(k)}_{0j}(\cdot)$ is continuous and its inverse, denoted by $H^{-(k)}_{0j}(\cdot)$, exists for all $j=1,\cdots,L_k$ and $k=1,2$. 
\end{theorem}
\begin{proof}
Note that the equality of joint sub-distribution functions, as per {\bf Definition 4}, implies equality of the unconditional $j$th sub-distribution functions $F^{(k)}_{j}(t_k;\boldsymbol{h^{(k)}_{0}},\sigma^{2}_{k})$ and $F^{(k)}_{j}(t_k;\boldsymbol{\tilde{h}^{(k)}_{0}},\tilde{\sigma}_k^2)$. Therefore, following the same technique as that used in the proof of {\bf Theorem 4.1} for the shared Gamma frailty model in the previous section, we have 
\begin{center}
    $\displaystyle{h^{(k)}_{0j}(t_k) = \tilde{h}^{(k)}_{0j}(t_k)}$ and 
    $\displaystyle{\sigma_{k} = \tilde{\sigma}_{k}},$ \\ 
\end{center}
for all $t_k>0,\ j=1,\cdots,L_k$ and $k = 1,2$. \\

Now, the equality of joint sub-distribution functions, as per {\bf Definition 4}, also implies equality of the unconditional joint survival functions $S(t_1,t_2;\boldsymbol{h^{(1)}_{0}},\boldsymbol{h^{(2)}_{0}},\boldsymbol{\theta})$ and $S(t_1,t_2;\boldsymbol{\tilde{h}^{(1)}_{0}},\boldsymbol{\tilde{h}^{(2)}_{0}},\boldsymbol{\tilde{\theta}})$. Note that the unconditional joint survival function 
can be obtained as 
\begin{align*}
S(t_1,t_2;\boldsymbol{h^{(1)}_{0}},\boldsymbol{h^{(2)}_{0}},\boldsymbol{\theta}) &= 
\int\limits_{0}^{\infty}\int\limits_{0}^{\infty}\int\limits_{0}^{\infty}\exp{\bigg[-\mu_{0}y_{0}\bigg(1+\sum_{k=1}^2 
\frac{H^{(k)}_{0}(t_k)}{\mu_k}\bigg) - \sum_{k=1}^2 
\mu_{k}y_{k}\bigg(1+\frac{H^{(k)}_{0}(t_k)}{\mu_k}\bigg)\bigg]} \\
&\qquad \times \prod_{l=0,1,2}\left(\frac{\mu^{\kappa_l}_l}{\Gamma(\kappa_l)}y^{\kappa_{l} - 1}_{l}dy_{l}\right) \\
&= \bigg(1+\sum_{k=1}^2 
\frac{H^{(k)}_{0}(t_k)}{\mu_k}\bigg)^{-\kappa_0}\times \prod_{k=1}^2 
\bigg(1+\frac{H^{(k)}_{0}(t_k)}{\mu_k}\bigg)^{-\kappa_k} \\
&= \Big(1+\sum_{k=1}^2 \sigma^{2}_{k}H^{(k)}_{0}(t_k)
\Big)^{-\frac{\rho}{\sigma_{1}\sigma_{2}}}\times \prod_{k=1}^2 
\Big(1+\sigma^{2}_{k}H^{(k)}_{0}(t_k)\Big)^{-\frac{1}{\sigma^{2}_k} + \frac{\rho}{\sigma_{1}\sigma_{2}}}, 
\end{align*}
for all $t_{1},t_{2} > 0$. Now, using the above expression in  the equality of $S(t_1,t_2;\boldsymbol{h^{(1)}_{0}},\boldsymbol{h^{(2)}_{0}},\boldsymbol{\theta})$ and $S(t_1,t_2;\boldsymbol{\tilde{h}^{(1)}_{0}},\boldsymbol{\tilde{h}^{(2)}_{0}},\boldsymbol{\tilde{\theta}})$, along with the already obtained equalities  
$\displaystyle{h^{(k)}_{0j}(t_k) = \tilde{h}^{(k)}_{0j}(t_k)}$, for all $t_{k}>0, j = 1,2,\cdots,L_{k}$, and 
$\displaystyle{\sigma_{k} = \tilde{\sigma}_{k}}$, for $k = 1,2$, we have 
\begin{multline*}
    \Big(1+\sum_{k=1}^2 \sigma^{2}_{k}H^{(k)}_{0}(t_k)\Big)^{-\frac{\rho}{\sigma_{1}\sigma_{2}}}\times \prod_{k=1}^2 
    \Big(1+\sigma^{2}_{k}H^{(k)}_{0}(t_k)\Big)^{-\frac{1}{\sigma^{2}_k} + \frac{\rho}{\sigma_{1}\sigma_{2}}} \\
\qquad =  \Big(1+\sum_{k=1}^2 \sigma^{2}_{k}H^{(k)}_{0}(t_k)\Big)^{-\frac{\tilde{\rho}}{\sigma_{1}\sigma_{2}}}\times \prod_{k=1}^2 
    \Big(1+\sigma^{2}_{k}H^{(k)}_{0}(t_k)\Big)^{-\frac{1}{\sigma^{2}_k} + \frac{\tilde{\rho}}{\sigma_{1}\sigma_{2}}}, 
\end{multline*}
for all $t_{1},t_{2} > 0$. This implies 
$$\Bigg[\frac{\Big(\sigma^{2}_{1}H^{(1)}_{0}(t_1) + 1\Big)\Big(\sigma^{2}_{2}H^{(2)}_{0}(t_2) + 1\Big)}{\Big(\sigma^{2}_{1}H^{(1)}_{0}(t_1) + \sigma^{2}_{2}H^{(2)}_{0}(t_2) + 1\Big)}\Bigg]^{\frac{\rho - \tilde{\rho}}{\sigma_{1}\sigma_{2}}} = 1,$$
for all $t_{1},t_{2} > 0$. Hence, $\rho = \tilde{\rho}$. \\
Therefore, the correlated Gamma frailty model $(\ref{NonparametrichazardcorrelatedGamma})$ is identifiable.
\end{proof}

\section{Shared cause-specific Gamma frailty model} \label{Nonparametric_hazard_shared_cause_specific_Gamma_frailty} 

In contrast with the individual level shared frailty and correlated frailty model of Sections 4 and 5, respectively, the frailty of the present section models individual frailty at the level of competing risks in the sense that the frailty variable depends on the cause of failure. We term this kind of frailty as cause-specific frailty. In this section, we consider shared cause-specific frailty in which two individuals in a pair share a common frailty for a particular cause of failure. As a result, this model requires same set of competing risks for the two individuals in a pair; that is, $L_1=L_2=L$, say. So, we have a common frailty variable $\epsilon_j$ for the two individuals for failure due to cause $j$, $j=1,\cdots,L$. While inducing dependence between the two individuals through the common $\epsilon_j$, for $j=1,\cdots,L$, it also allows different frailty for different causes. In the example of multi-organ failure, while possible dependence between the organs is captured by the common frailty variables, the $\epsilon_j$'s, the nature of frailty depends on the cause of failure. \\

Formally, the shared cause-specific frailty model is given by 
\begin{equation}\label{shared_cause_specificGamma}
   \lambda^{(k)}_{j}(t_k | \boldsymbol{\epsilon}) = h^{(k)}_{0j}(t_k) \epsilon_{j},
\end{equation}
for all $t_k>0,\ j = 1,\cdots,L$ and $k=1,2$, where  $\boldsymbol{\epsilon} = (\epsilon_1,\epsilon_2,\cdots,\epsilon_L)$.  We assume that the cause-specific frailty variables $\epsilon_{j}$'s are independent with $\epsilon_{j}$ following Gamma$(\frac{1}{\sigma_j^2},\frac{1}{\sigma_j^2})$ distribution with $\sigma_j>0$ having density same as that in (\ref{gammadensity}), so that $E(\epsilon_{j})=1$, for $j=1,\cdots,L$. Therefore, the frailty parameter space $\boldsymbol{\Theta}$ is given by 
$$\boldsymbol{\Theta}=\big\{\boldsymbol{\theta}=(\sigma_1,\cdots,\sigma_L): \sigma_{j} > 0, j=1,\cdots,L\big\}.$$\label{Theta_3}. 
Then, as before, the survival function of the $k$th individual, conditional on $\boldsymbol{\epsilon}$, is
$$S^{(k)}(t_{k};\boldsymbol{h^{(k)}_{0}}|\boldsymbol{\epsilon}) 
    = \exp{\bigg[- \sum\limits_{j=1}^{L}H^{(k)}_{0j}(t_k)\epsilon_{j}\bigg]},$$
and the $j$th sub-distribution function of the $k$th individual, conditional on $\boldsymbol{\epsilon}$, is 
$$F^{(k)}_{j}(t_k;\boldsymbol{h^{(k)}_{0}}| \boldsymbol{\epsilon})= \int\limits_{0}^{t_k}h^{(k)}_{0j}(u_k)\epsilon_{j}\exp{\bigg[- \sum\limits_{j'=1}^{L}H^{(k)}_{0j'}(u_k)\epsilon_{j'}\bigg]}du_k,
$$
for all $t_{k} > 0,\ j = 1,\cdots,L$ and $k = 1,2$. 
The unconditional joint sub-distribution function $F_{j_1j_2}(t_1,t_2;\boldsymbol{h^{(1)}_{0}},\boldsymbol{h^{(2)}_{0}};\boldsymbol{\theta})$ is given by 
$$\int\limits_{0}^{\infty}\cdots\int\limits_{0}^{\infty}\Bigg[\int\limits_{0}^{t_1}\int\limits_{0}^{t_2} \left(
h^{(k)}_{0j_k}(u_k)\epsilon_{j_k}\right) \times 
\exp{\bigg(-\sum\limits_{j=1}^{L}\epsilon_{j}\sum_{k=1}^2H^{(k)}_{0j}(u_k)\bigg)} du_2du_1\Bigg]g(\boldsymbol{\epsilon};\boldsymbol{\theta})d \boldsymbol{\epsilon}, $$
for all $t_k>0,\ j_k=1,\cdots,L$ and $k=1,2$, 
where $g(\boldsymbol{\epsilon};\boldsymbol{\theta})$ denotes the product density of $L$ independent Gamma densities for $\epsilon_1,\cdots,\epsilon_L$, as described above. \\ 
The unconditional sub-distribution function $F^{(k)}_{j}(t_k;\boldsymbol{h^{(k)}_{0}}, \boldsymbol{\theta})$ is given by 
\begin{align}\label{Subdistri_shared_cause_specific}
    &= \int\limits_{0}^{t_k}h^{(k)}_{0j}(u_k)\Bigg[\int\limits_{0}^{\infty}\cdots\int\limits_{0}^{\infty} \epsilon_j \exp{\bigg(- \sum\limits_{j'=1}^{L}H^{(k)}_{0j'}(u_k)\epsilon_{j'}\bigg)}\prod\limits_{j'=1}^{L}\frac{e^{-\frac{\epsilon_{j'}}{\sigma_{j'}^2}}\epsilon_{j'}^{\frac{1}{\sigma_{j'}^2} - 1}}{\sigma_{j'}^{\frac{2}{\sigma_{j'}^2}}\Gamma(\frac{1}{\sigma_{j'}^2})} d\epsilon_{j'}\Bigg]du_k \nonumber \\
    &= \int\limits_{0}^{t_k}h^{(k)}_{0j}(u_k) \Bigg[ \Bigg(
    \frac{1}{\sigma_j^{\frac{2}{\sigma_{j}^2}} \Gamma(\frac{1}{\sigma_j^2})}\int\limits_{0}^{\infty} \exp{\bigg(-
    \epsilon_j
    \Big(H^{(k)}_{0j}(u_k) + \frac{1}{\sigma_j^2}\Big)\bigg)} 
    \epsilon_j^{\frac{1}{\sigma_j^2}}d\epsilon_j \Bigg) \times \nonumber \\
   & \qquad\prod\limits_{j'(\ne j) = 1}^{L}\Bigg(\frac{1}{\sigma_{j'}^{\frac{2}{\sigma_{j'}^2}}\Gamma(\frac{1}{\sigma_{j'}^2})}\int\limits_{0}^{\infty}\exp\bigg(- \epsilon_{j'}\Big(H^{(k)}_{0j'}(u) + \frac{1}{\sigma_{j'}^2}\Big)\bigg)\epsilon_{j'}^{\frac{1}{\sigma_{j'}^2} - 1} d\epsilon_{j'} \Bigg) \Bigg]du_k \nonumber \\
  &= \int\limits_{0}^{t_k}h^{(k)}_{0j}(u_k) \Bigg[ \Bigg( 
    \frac{1}{\sigma_j^{\frac{2}{\sigma_{j}^2}} \Gamma(\frac{1}{\sigma_j^2})} \frac{\Gamma(\frac{1}{\sigma_j^2}+1)} {\Big(H^{(k)}_{0j}(u_k) + \frac{1}{\sigma_j^2}\Big)^{\frac{1}{\sigma^{2}_j}+1}}\Bigg)  \times  
    \prod\limits_{j'(\ne j) = 1}^{L}\Bigg(\frac{1}{\sigma_{j'}^{\frac{2}{\sigma_{j'}^2}}\Gamma(\frac{1}{\sigma_{j'}^2})} \times 
    \frac{\Gamma(\frac{1}{\sigma_{j'}^2})} {\Big(H^{(k)}_{0j'}(u_k) + \frac{1}{\sigma^{2}_{j'}}\Big)^{\frac{1}{\sigma^{2}_{j'}}}}\Bigg) \Bigg]du_k \nonumber \\
    &= \int\limits_{0}^{t_k}h^{(k)}_{0j}(u_k) \Bigg[ \Bigg( 
   \frac{1}{\Big(1+\sigma_j^2H^{(k)}_{0j}(u_k)\Big)^{\frac{1}{\sigma_j^2}+1}}\Bigg)  \times  \prod\limits_{j'(\ne j) = 1}^{L}\Bigg(\frac{1} {\Big(1+\sigma_{j'}^2H^{(k)}_{0j'}(u_k)\Big)^{\frac{1}{\sigma^{2}_{j'}}}}\Bigg) \Bigg]du_k \nonumber \\
 &= \int\limits_{0}^{t_k}h^{(k)}_{0j}(u_k) \Bigg[ 
   \frac{1}{1+\sigma_j^2H^{(k)}_{0j}(u_k)} \times  \prod\limits_{j' = 1}^{L}\Bigg(\frac{1} {\Big(1+\sigma^{2}_{j'}H^{(k)}_{0j'}(u_k)\Big)^{\frac{1}{\sigma^{2}_{j'}}}}\Bigg) \Bigg]du_k, 
\end{align}
for all $t_k>0,\ j=1,\cdots,L$ and $k=1,2$. 

\begin{definition}
The shared cause-specific Gamma frailty model (\ref{shared_cause_specificGamma}) is identifiable within the family $\mathcal{H}_1\times\mathcal{H}_2\times\boldsymbol{\Theta}$ if, for some $\boldsymbol{h^{(k)}_{0}},\boldsymbol{\tilde{h}^{(k)}_{0}} \in \mathcal{H}_{k}$, for $k=1,2$, and $\boldsymbol{\theta},\boldsymbol{\tilde{\theta}} \in \boldsymbol{\Theta}$, where $\boldsymbol{\theta}=(\sigma_1,\cdots,\sigma_L),\ \boldsymbol{\tilde{\theta}}=(\tilde{\sigma}_1,\cdots,\tilde{\sigma}_L),$ the equality of the joint sub-distribution functions 
$$F_{j_1j_2}(t_1,t_2;\boldsymbol{h^{(1)}_{0}},\boldsymbol{h^{(2)}_{0}};\boldsymbol{\theta})= F_{j_1j_2}(t_1,t_2;\boldsymbol{\tilde{h}^{(1)}_{0}},\boldsymbol{\tilde{h}^{(2)}_{0}};\boldsymbol{\tilde{\theta}}),$$
for all $t_k>0,\ j_k=1,\cdots,L$ and $k=1,2$, implies 
$$h^{(k)}_{0j}(x) = \tilde{h}^{(k)}_{0j}(x) \, \, \text{for all} \, \, x>0,\ \ \sigma_{j} = \tilde{\sigma}_{j}, \, \, \text{for} \, \ j= 1,\cdots,L\, \, \text{and} \, \, k = 1,2.$$
\end{definition}
\begin{theorem}\label{Theorem 4.3}
The shared cause-specific Gamma frailty model (\ref{shared_cause_specificGamma}) with 
non-parametric baseline cause-specific hazard functions is identifiable within the family $\mathcal{H}_{1}\times\mathcal{H}_2\times\boldsymbol{\Theta}$, where  $\boldsymbol{\Theta}=\big\{\boldsymbol{\theta}=(\sigma_1,\cdots,\sigma_L): \sigma_{1} > 0, \cdots,\sigma_{L} > 0\big\}$,  provided $H^{(k)}_{0j}(\cdot)$ is continuous and its inverse, denoted by $H^{-(k)}_{0j}(\cdot)$, exists for all $j=1,\cdots,L$ and $k=1,2$.

\end{theorem}

\begin{proof} 

As before, the equality of joint sub-distribution functions, as per {\bf Definition 4}, implies equality of the unconditional $j$th sub-distribution functions $F^{(k)}_{j}(t_k;\boldsymbol{h^{(k)}_{0}},\boldsymbol{\theta})$ and $F^{(k)}_{j}(t_k;\boldsymbol{\tilde{h}^{(k)}_{0}},\tilde{\boldsymbol{\theta}})$, as given by (\ref{Subdistri_shared_cause_specific}), for all $t_k>0,\ j=1,\cdots,L$ and $k=1,2$. Hence, we have equality of the corresponding sub-density functions, as given by  
\begin{multline}\label{shared_cause_specific_identifiabilityEq_2}
  h^{(k)}_{0j}(t_k) \big\{1 + \sigma_{j}^2H^{(k)}_{0j}(t_k)\big\}^{-1} 
  \prod\limits_{j'=1}^{L}\big\{1 + \sigma_{j'}^2 H^{(k)}_{0j'}(t_k)\big\}^{-\frac{1}{\sigma_{j'}^2}}
 \\ =\tilde{h}^{(k)}_{0j}(t_k) \big\{1 + \tilde{\sigma}_j^2\tilde{H}^{(k)}_{0j}(t_k)\big\}^{-1} 
 \prod\limits_{j'=1}^{L}\big\{1 + \tilde{\sigma}_{j'}^2 \tilde{H}^{(k)}_{0j'}(t_k)\big\}^{-\frac{1}{\tilde{\sigma}_{j'}^2}},
\end{multline}
for all $t_k>0,\ j=1,\cdots,L$ and $k=1,2$. This implies 
$$
 \frac{\frac{d}{dt_k}\Big[\frac{1}{\sigma_{j}^2}\log(1 + \sigma_{j}^2H^{(k)}_{0j}(t_k))\Big]}{\frac{d}{dt_k}\Big[\frac{1}{\tilde{\sigma}_{j}^2}\log(1 + \tilde{\sigma}_{j}^2\tilde{H}^{(k)}_{0j}(t_k))\Big]} = 
 \frac{\big\{1 + \tilde{\sigma}_{j}^2 \tilde{H}^{(k)}_{0j}(t_k)\big\}^{-\frac{1}{\tilde{\sigma}_{j}^2}} 
 \prod\limits_{j'(\ne j)=1}^{L}\big\{1 + \tilde{\sigma}_{j'}^2 \tilde{H}^{(k)}_{0j'}(t_k)\big\}^{-\frac{1}{\tilde{\sigma}_{j'}^2}}}
 {\big\{1 + \sigma_{j}^2 H^{(k)}_{0j}(t_k)\big\}^{-\frac{1}{\sigma_{j}^2}}
 \prod\limits_{j'(\ne j)=1}^{L}\big\{1 + \sigma_{j'}^2 H^{(k)}_{0j'}(t_k)\big\}^{-\frac{1}{\sigma_{j'}^2}}},$$
 for all $t_k>0,\ j=1,\cdots,L$ and $k=1,2$. 
 Let us now write $y=\log[1 + \sigma_{j}^2 H^{(k)}_{0j}(t_k)]$\\ 
 and $z=\log[1 + \tilde{\sigma}_{j}^2 \tilde{H}^{(k)}_{0j}(t_k)]$. Then, from the above equality, we have 
 $$\frac{\tilde{\sigma}_j^2}{\sigma_j^2}\times \frac{dy}{dz}= 
 \frac{e^{-\frac{z}{\tilde{\sigma}_j^2}} \prod\limits_{j' (\ne j)=1}^{L} \big[1 + \tilde{\sigma}_{j'}^2 \tilde{H}^{(k)}_{0j'}\left(\tilde{H}^{-(k)}_{0j}\big(\frac{e^z-1}{\tilde{\sigma}_j^2}\big)\right)\big]^{-\frac{1}{\tilde{\sigma}_{j'}^2}}}
 {e^{-\frac{y}{\sigma_j^2}} \prod\limits_{j' (\ne j)=1}^{L} \big[1 + \sigma_{j'}^2 H^{(k)}_{0j'}\left(H^{-(k)}_{0j}\big(\frac{e^y-1}{\sigma_j^2}\big)\right)\big]^{-\frac{1}{\sigma_{j'}^2}}}. $$

This is equivalent to, after integrating from 0 to some $w>0$, 
$$
\frac{\int_0^w e^{-\frac{y}{\sigma_j^2}} \prod\limits_{j' (\ne j)=1}^{L} \big[1 + \sigma_{j'}^2 H^{(k)}_{0j'}\left(H^{-(k)}_{0j}\big(\frac{e^y-1}{\sigma_j^2}\big)\right)\big]^{-\frac{1}{\sigma_{j'}^2}}dy}
{\int_0^w e^{-\frac{z}{\tilde{\sigma}_j^2}} \prod\limits_{j' (\ne j)=1}^{L} \big[1 + \tilde{\sigma}_{j'}^2 \tilde{H}^{(k)}_{0j'}\left(\tilde{H}^{-(k)}_{0j}\big(\frac{e^z-1}{\tilde{\sigma}_j^2}\big)\right)\big]^{-\frac{1}{\tilde{\sigma}_{j'}^2}}dz}= \frac{\sigma_j^2}{\tilde{\sigma}_j^2}.$$

Now, taking limit as $w \to 0+$ and applyling L'opital rule, we get 
$$\lim\limits_{w \to 0+}\frac{
e^{-\frac{w}{\sigma_j^2}} \prod\limits_{j' (\ne j)=1}^{L} \big[1 + \sigma_{j'}^2 H^{(k)}_{0j'}\left(H^{-(k)}_{0j}\big(\frac{e^w-1}{\sigma_j^2}\big)\right)\big]^{-\frac{1}{\sigma_{j'}^2}}}
{e^{-\frac{w}{\tilde{\sigma}_j^2}} \prod\limits_{j' (\ne j)=1}^{L} \big[1 + \tilde{\sigma}_{j'}^2 \tilde{H}^{(k)}_{0j'}\left(\tilde{H}^{-(k)}_{0j}\big(\frac{e^w-1}{\tilde{\sigma}_j^2}\big)\right)\big]^{-\frac{1}{\tilde{\sigma}_{j'}^2}}}= \frac{\sigma_j^2}{\tilde{\sigma}_j^2}.$$
 Note that, as in Section 4, the composite functions appearing in the above equality are continuous functions. Hence, for $j'\ne j$, 
 $$\lim\limits_{w \to 0+}
H^{(k)}_{0j'}\left(H^{-(k)}_{0j}\big(\frac{e^w-1}{\sigma_j^2}\big)\right)= \lim\limits_{w \to 0+}
\tilde{H}^{(k)}_{0j'}\left(\tilde{H}^{-(k)}_{0j}\big(\frac{e^w-1}{\tilde{\sigma}_j^2}\big)\right)=0.$$
Hence, since $\lim\limits_{w \to 0+}e^{-\frac{w}{\sigma_j^2}}=  
\lim\limits_{w \to 0+}e^{-\frac{w}{\tilde{\sigma}_j^2}} = 1$, 
we get $\frac{\sigma_j^2}{\tilde{\sigma}_j^2}=1$, or $\sigma_j^2=\tilde{\sigma}_j^2$, for $j=1,\cdots,L$.\\
Using this equality in (\ref{shared_cause_specific_identifiabilityEq_2}) and then summing over $j=1,\cdots,L$, we get 
\begin{multline*}
    \left(\prod\limits_{j'=1}^{L}\big[1 + \sigma_{j'}^2 H^{(k)}_{0j'}(t_k)\big]^{-\frac{1}{\sigma_{j'}^2}}\right)\times\left(
    \sum\limits_{j=1}^{L}\frac{h^{(k)}_{0j}(t_k)}{1 + \sigma_{j}^2H^{(k)}_{0j}(t_k)}\right)\\
    = \left( \prod\limits_{j'=1}^{L}\big[1 + \sigma_{j'}^2 \tilde{H}^{(k)}_{0j'}(t_k)\big]^{-\frac{1}{\sigma_{j'}^2}}\right)\times\left(
\sum\limits_{j=1}^{L}\frac{\tilde{h}^{(k)}_{0j}(t_k)}{1 + \sigma_{j}^2\tilde{H}^{(k)}_{0j}(t_k)}\right), 
\end{multline*}
for all $t_k>0$ and $k=1,2$. Writing 
$$c(t_k)=\log\prod\limits_{j'=1}^{L}\big[1 + \sigma_{j'}^2 H^{(k)}_{0j'}(t_k)\big]^{\frac{1}{\sigma_{j'}^2}}\quad\mbox{ and }\quad \tilde{c}(t_k)=\log\prod\limits_{j'=1}^{L}\big[1 + \sigma_{j'}^2 \tilde{H}^{(k)}_{0j'}(t_k)\big]^{\frac{1}{\sigma_{j'}^2}},  $$
the above equality reduces to 
$$e^{-c(t_k)}\frac{d}{dt_k}c(t_k)=e^{-\tilde{c}(t_k)}\frac{d}{dt_k}\tilde{c}(t_k),$$
for all $t_k>0$.  By integrating both sides, we get $e^{-c(t_k)}=e^{-\tilde{c}(t_k)}$, since $\lim\limits_{t_{k} \to 0+}c(t_k) =  \lim\limits_{t_{k} \to 0+}\tilde{c}(t_k) = 0$. This implies  $c(t_k)=\tilde{c}(t_k)$ for all $t_k>0$. That is, 
$$\prod\limits_{j'=1}^{L}\big[1 + \sigma_{j'}^2 H^{(k)}_{0j'}(t_k)\big]^{\frac{1}{\sigma_{j'}^2}}=\prod\limits_{j'=1}^{L}\big[1 + \sigma_{j'}^2 \tilde{H}^{(k)}_{0j'}(t_k)\big]^{\frac{1}{\sigma_{j'}^2}}.$$
Using this equality in (\ref{shared_cause_specific_identifiabilityEq_2}) and then integrating, we get 
$$\int\limits_{0}^{t_k}h^{(k)}_{0j}(u_k)
\big[1 + \sigma_{j}^2H^{(k)}_{0j}(u_k)\big]^{-1} du_k  =  \int\limits_{0}^{t_k}\tilde{h}^{(k)}_{0j}(u_k)
\big[1 + \sigma_{j}^2\tilde{H}^{(k)}_{0j}(u_k)\big]^{-1} du_k.$$
This gives $H^{(k)}_{0j}(t_k)=\tilde{H}^{(k)}_{0j}(t_k)$, or 
$h^{(k)}_{0j}(t_k)=\tilde{h}^{(k)}_{0j}(t_k)$, for all $t_k>0,\  j=1,\cdots,L$ and $k=1,2$. Hence, the theorem is proved.  
\end{proof}

\section{Correlated cause-specific Gamma frailty model} \label{Nonparametric_hazard_correlated_cause_specific_Gamma_frailty}

As in the previous section, the frailty of the present section models individual frailty at the level of competing risks; but, instead of being shared between the two individuals, there are two different correlated frailty variables for each cause of failure. That is, for a particular cause of failure, say $j$th, there are two correlated frailty variables $\epsilon_j^{(k)},\ k=1,2$, thus allowing for different frailty for different individuals and different causes. Since the frailties are cause-specific, we need to have same set of competing risks for both the individuals, that is $L_1=L_2=L$, say. 
In the example of multi-organ failure, while possible dependence between the organs is captured by the correlated frailty variables, the $\epsilon_j^{(k)},\ k=1,2$, the nature of frailty depends on the cause of failure. \\

Formally, the correlated cause-specific frailty model is given by
\begin{equation}\label{correlated_cause_specific_Gamma_Model}
   \lambda^{(k)}_{j}(t_k | \boldsymbol{\epsilon^{(k)}}) = h^{(k)}_{0j}(t_k)\epsilon^{(k)}_{j},
\end{equation}
where $\boldsymbol{\epsilon^{(k)}}= (\epsilon^{(k)}_{1},\cdots,\epsilon^{(k)}_{L})$, for all $t_k>0,\ j= 1,\cdots,L$ and $k=1,2$. We assume that, while $(\epsilon_j^{(1)}, \epsilon_j^{(2)})$ follows a correlated Gamma frailty distribution for a particular $j$ (See Section 5), the $L$ pairs of frailty variables $(\epsilon_j^{(1)}, \epsilon_j^{(2)})$, for $ j=1,\cdots,L,$ are independent. As in Section 5, for a particular $j$, write 
$\epsilon^{(k)}_j = \frac{\mu_{0j}}{\mu_{kj}}Y_{0j} + Y_{kj}$, for $k=1,2$, where $Y_{0j} \sim \text{Gamma}(\kappa_{0j},\mu_{0j}),Y_{1j} \sim \text{Gamma}(\kappa_{1j},\mu_{1j})$ and $Y_{2j} \sim \text{Gamma}(\kappa_{2j},\mu_{2j})$, with  $\mu_{0j},\mu_{1j},\mu_{2j},\kappa_{0j},\kappa_{1j},\kappa_{2j} > 0$ satisfying $\mu_{kj} = \kappa_{0j} + \kappa_{kj}$, for $k=1,2$, and the random variables $Y_{0j},Y_{1j},Y_{2j}$ are independent. For a particular $j$, we have three independent parameters given by 
$$\sigma_{1j}=\frac{1}{\sqrt{\kappa_{0j} + \kappa_{1j}}},\ \sigma_{2j}=\frac{1}{\sqrt{\kappa_{0j} + \kappa_{2j}}},\ \rho_j=\frac{\kappa_{0j}}{\sqrt{\kappa_{0j} + \kappa_{1j}}\sqrt{\kappa_{0j} + \kappa_{2j}}},$$
so that $\sigma_{kj}>0$ is the standard deviation of the frailty variable $\epsilon^{(k)}_j$, for $k=1,2$ and $\rho_j$ is the  correlation coefficient between $\epsilon^{(1)}_j$ and $\epsilon^{(2)}_j$ satisfying $0 < \rho_j < \min{(\frac{\sigma_{1j}}{\sigma_{2j}},\frac{\sigma_{2j}}{\sigma_{1j}})}$, for $j=1,\cdots,L$. 
Therefore, as in Section 5, the frailty parameter space for correlated cause-specific Gamma frailty model is given by 
$$\boldsymbol{\Theta}=\big\{\boldsymbol{\theta}=(\sigma_{1j},\sigma_{2j},\rho_j): \sigma_{1j} > 0, \sigma_{2j} > 0,\  0 < \rho_j < \min{(\frac{\sigma_{1j}}{\sigma_{2j}},\frac{\sigma_{2j}}{\sigma_{1j}})},\ j=1,\cdots,L\big\}.$$\label{Theta_5}.

The survival function of the $k$th individual, conditional on the frailty vector $\boldsymbol{\epsilon^{(k)}}$, is 
$$S^{(k)}(t_{k};\boldsymbol{h^{(k)}_0}| \boldsymbol{\epsilon^{(k)}}) 
    = \exp{\bigg[-\sum\limits_{j = 1}^{L}H^{(k)}_{0j}(t_k)\epsilon^{(k)}_{j}\bigg]},$$
and the $j$th sub-distribution function of the $k$th individual, conditional on frailty vector $\boldsymbol{\epsilon^{(k)}}$, is 
$$ F^{(k)}_{j}(t_{k};\boldsymbol{h^{(k)}_0} | \boldsymbol{\epsilon^{(k)}}) = \int\limits_{0}^{t_k}h^{(k)}_{0j}(u_k)\epsilon^{(k)}_{j}\exp{\bigg[-\sum\limits_{j' = 1}^{L}H^{(k)}_{0j'}(u_k)\epsilon^{(k)}_{j'}\bigg]}du_k, $$
for all $t_{k} > 0,\ j = 1,\cdots,L$ and $k = 1,2$. The unconditional joint sub-distribution $F_{j_1j_2}(t_1,t_2;\boldsymbol{h^{(1)}_0},\boldsymbol{h^{(2)}_0}, \boldsymbol{\theta})$ is 
$$\int\limits_{0}^{\infty}\cdots\int\limits_{0}^{\infty}\Bigg[\int\limits_{0}^{t_1}\int\limits_{0}^{t_2} \prod_{k=1}^2\left(
\epsilon^{(k)}_{j_k} h^{(k)}_{0j_k}(u_k)\right)\times 
\exp{\bigg(- \sum_{k=1}^2 
\sum\limits_{j' = 1}^{L}H^{(k)}_{0j'}(u_k) \epsilon^{(k)}_{j'}\bigg)} du_2 du_1\Bigg]
g(\boldsymbol{\epsilon^{(1)}},\boldsymbol{\epsilon^{(2)}}; \boldsymbol{\theta}) d\boldsymbol{\epsilon^{(1)}}d\boldsymbol{\epsilon^{(2)}}, $$
for all $t_{k} > 0,\ j_{k} = 1,\cdots,L$ and $k = 1,2$, where 
$g(\boldsymbol{\epsilon^{(1)}},\boldsymbol{\epsilon^{(2)}}; \boldsymbol{\theta})$ denotes the joint density of $(\boldsymbol{\epsilon^{(1)}},\boldsymbol{\epsilon^{(2)}})$ following a correlated cause-specific Gamma distribution, as described above. \\

As in Section 5, from the expression of the $j$th sub-distribution function of the $k$th individual, conditional on $\boldsymbol{\epsilon^{(k)}}$, or equivalently $\{(Y_{0j}, Y_{kj}),\ j=1,\cdots,L\}$, we get the corresponding  unconditional $j$th sub-distribution function, for $k=1,2$, as 
\begin{align*}
    F^{(k)}_{j}(t_k;\boldsymbol{h^{(k)}_0};\boldsymbol{\theta}) &=    \int\limits_{0}^{t_k}\int\limits_{0}^{\infty}\int\limits_{0}^{\infty}\cdots\int\limits_{0}^{\infty}h^{(k)}_{0j}(u_k)\bigg(\frac{\mu_{0j}}{\mu_{kj}}y_{0j} + y_{kj}\bigg)\exp{\bigg[-\sum\limits_{j'=1}^{L}H^{(k)}_{0j'}(u_k)\big(\frac{\mu_{0j'}}{\mu_{kj'}}y_{0j'} + y_{kj'}\big)\bigg]}\times \\
    &\hspace{2.0in} \prod\limits_{l=0,k}\prod\limits_{j'=1}^{L} \left(\frac{\mu_{lj'}^{\kappa_{lj'}}}{\Gamma(\kappa_{lj'})}e^{-\mu_{lj'}y_{lj'}} y_{lj'}^{\kappa_{lj'} - 1} dy_{lj'}\right)du_k \\
    &= \int\limits_{0}^{t_k}h^{(k)}_{0j}(u_k)\big(I_{1j}(u_k) + I_{2j}(u_k)\big)du_k \qquad\text{(say)},
\end{align*}
where $I_{1j}(u_k)$ is given by 
\begin{align*}
     &\left( \frac{\mu_{0j}}{\mu_{kj}}\prod\limits_{l=0,k}\prod\limits_{j'=1}^{L} \frac{\mu_{lj'}^{\kappa_{lj'}}}{\Gamma(\kappa_{lj'})}\right)\int\limits_{0}^{\infty}\cdots\int\limits_{0}^{\infty}\exp{\bigg[-\sum\limits_{j'=1}^{L}\Big(1+\frac{H^{(k)}_{0j'}(u_k)}{\mu_{kj'}}\Big)\big(\sum_{l=0,k} y_{lj'}\mu_{lj'}\big)\bigg]} 
     y_{0j}\prod\limits_{j'=1}^{L}\prod_{l=0,k}\left(
     y^{\kappa_{lj'} - 1}_{lj'}dy_{lj'}\right) \\ 
    &\qquad = \left(\frac{\mu_{0j}}{\mu_{kj}} \frac{\mu_{0j}^{\kappa_{0j}}}{\Gamma(\kappa_{0j})}\int\limits_{0}^{\infty}\exp{\bigg[-\Big(1+\frac{H^{(k)}_{0j}(u_k)}{\mu_{kj}}\Big) y_{0j}\mu_{0j}\bigg]}y^{\kappa_{0j} }_{0j}dy_{0j}\right) \\ 
    &\qquad\times \left(
    \prod\limits_{j' (\ne j)=1}^{L}\frac{\mu_{0j'}^{\kappa_{0j'}}}{\Gamma(\kappa_{0j'})}\int\limits_{0}^{\infty}\exp\bigg[-\Big(1+\frac{H^{(k)}_{0j'}(u_k)}{\mu_{kj'}}\Big) y_{0j'}\mu_{0j'}\bigg]y^{\kappa_{0j'} - 1}_{0j'}dy_{0j'}\right) \\
    &\qquad\times \left(\prod\limits_{j'=1}^{L} \frac{\mu_{kj'}^{\kappa_{kj'}}}{\Gamma(\kappa_{kj'})}\int\limits_{0}^{\infty}\exp\bigg[-\Big(1+\frac{H^{(k)}_{0j'}(u_k)}{\mu_{kj'}} \Big) y_{kj'}\mu_{kj'}\bigg]y^{\kappa_{kj'} - 1}_{kj'}dy_{kj'}\right) \\
    &\qquad = \left(\frac{\mu_{0j}}{\mu_{kj}} 
    \frac{\kappa_{0j}}{\mu_{0j}}\Big(1+\frac{H^{(k)}_{0j}(u_k)}{\mu_{kj}}\Big)^{-\kappa_{0j} - 1}\right)\times\left(
    \prod\limits_{j' (\ne j)=1}^{L}\Big(1+\frac{H^{(k)}_{0j'}(u_k)}{\mu_{kj'}}\Big)^{-\kappa_{0j'}}\right)\times\left(
     \prod\limits_{j'=1}^{L}\Big(1+\frac{H^{(k)}_{0j'}(u_k)}{\mu_{kj'}}\Big)^{-\kappa_{kj'}}\right) \\
    &\qquad = \left(\frac{\kappa_{0j}}{\mu_{kj}}\prod\limits_{j'=1}^{L}\Big(1+\frac{H^{(k)}_{0j'}(u_k)}{\mu_{kj'}}\Big)^{-\kappa_{0j'}-\kappa_{kj'}}\right)\times  \Big(1+\frac{H^{(k)}_{0j}(u_k)}{\mu_{kj}}\Big)^{-1} \\
    &\qquad = \frac{\kappa_{0j}}{\mu_{kj}} \Big(1+\frac{H^{(k)}_{0j}(u_k)}{\mu_{kj}}\Big)^{-1} \prod\limits_{j'=1}^{L}\Big(1+\frac{H^{(k)}_{0j'}(u_k)}{\mu_{kj'}}\Big)^{-\mu_{kj'}}
\end{align*}
and $I_{2j}(u_k)$ is given by 
\begin{align*}
  &\left( \prod\limits_{l=0,k}\prod\limits_{j'=1}^{L} \frac{\mu_{lj'}^{\kappa_{lj'}}}{\Gamma(\kappa_{lj'})}\right) \int\limits_{0}^{\infty}\cdots\int\limits_{0}^{\infty}\exp{\bigg[-\sum\limits_{j'=1}^{L}\Big(1+\frac{H^{(k)}_{0j'}(u_k)}{\mu_{kj'}}\Big)\sum_{l=0,k}y_{lj'}\mu_{lj'}\bigg]} y_{kj} 
  \prod_{l=0,k}\prod\limits_{j'=1}^{L}y^{\kappa_{lj'} - 1}_{lj'}dy_{lj'} \\   
  &\qquad =\left(\frac{\mu_{kj}^{\kappa_{kj}}}{\Gamma(\kappa_{kj})}\int\limits_{0}^{\infty}\exp{\bigg[-\Big(1+\frac{H^{(k)}_{0j}(u_k)}{\mu_{kj}}\Big) y_{kj}\mu_{kj}\bigg]}y^{\kappa_{kj} }_{kj}dy_{kj}\right) \\
  &\qquad\times\left(\prod\limits_{j' (\ne j)=1}^{L} \frac{\mu_{kj'}^{\kappa_{kj'}}}{\Gamma(\kappa_{kj'})} \int\limits_{0}^{\infty}\exp{\bigg[-\Big(1+\frac{H^{(k)}_{0j'}(u_k)}{\mu_{kj'}}\Big)y_{kj'}\mu_{kj'}\bigg]}y^{\kappa_{kj'} - 1}_{kj'}dy_{kj'}\right) \\
  &\qquad\times\left(\prod\limits_{j'=1}^{L} \frac{\mu_{0j'}^{\kappa_{0j'}}}{\Gamma(\kappa_{0j'})} \int\limits_{0}^{\infty}\exp{\bigg[-\Big(1+\frac{H^{(k)}_{0j'}(u_k)}{\mu_{kj'}}\Big)y_{0j'}\mu_{0j'}\bigg]}y^{\kappa_{0j'} - 1}_{0j'}dy_{0j'}\right) \\
  &\qquad =\left( \frac{\kappa_{kj}}{\mu_{kj}}\Big(1+\frac{H^{(k)}_{0j}(u_k)}{\mu_{kj}}\Big)^{-\kappa_{kj}-1}\right)\times\left(
  \prod\limits_{j' (\ne j)=1}^{L}\Big(1+\frac{H^{(k)}_{0j'}(u_k)}{\mu_{kj'}}\Big)^{-\kappa_{kj'}} \right)\times\left( \prod\limits_{j'=1}^{L}\Big(1+\frac{H^{(k)}_{0j'}(u_k)}{\mu_{kj'}}\Big)^{-\kappa_{0j'}}\right) \\ 
  &\qquad = \frac{\kappa_{kj}}{\mu_{kj}} \Big(1+\frac{H^{(k)}_{0j}(u_k)}{\mu_{kj}}\Big)^{-1} \prod\limits_{j'=1}^{L}\Big(1+\frac{H^{(k)}_{0j'}(u_k)}{\mu_{kj'}}\Big)^{-\kappa_{0j'}-\kappa_{kj'}} \\
  &\qquad= \frac{\kappa_{kj}}{\mu_{kj}} \Big(1+\frac{H^{(k)}_{0j}(u_k)}{\mu_{kj}}\Big)^{-1} \prod\limits_{j'=1}^{L}\Big(1+\frac{H^{(k)}_{0j'}(u_k)}{\mu_{kj'}}\Big)^{-\mu_{kj'}}
\end{align*}
for all $u_k > 0$ and $k=1,2$. Therefore, we have 
\begin{align*}
   F^{(k)}_{j}(t_k;\boldsymbol{h^{(k)}_0},\boldsymbol{\theta}) &= \int\limits_{0}^{t_k}h^{(k)}_{0j}(u_k)\bigg[\frac{\kappa_{0j}+\kappa_{kj}}{\mu_{kj}}\bigg(1+\frac{H^{(k)}_{0j}(u_k)}{\mu_{kj}}\bigg)^{-1} \times\prod\limits_{j'=1}^{L}\Big(1+\frac{H^{(k)}_{0j'}(u_k)}{\mu_{kj'}}\Big)^{-\mu_{kj'}}\bigg]du_k \\
   &= \int\limits_{0}^{t_k}h^{(k)}_{0j}(u_k)\Bigg[\frac{1}{1 + \sigma^{2}_{kj}H^{(k)}_{0j}(u_k)}\prod\limits_{j'=1}^{L}\frac{1}{\Big(1 + \sigma^{2}_{kj'}H^{(k)}_{0j'}(u_k)\Big)^{\frac{1}{\sigma^{2}_{kj'}}}}\Bigg]du_k
\end{align*}
for all $t_{k} > 0,\ j= 1,\cdots,L$ and $k = 1,2$. 
Note that this unconditional $j$th sub-distribution function of the $k$th individual has the same form as that in (\ref{Subdistri_shared_cause_specific}) for the shared cause-specific Gamma frailty model of the previous section. 
\newpage
\begin{definition}
The correlated cause-specific Gamma frailty model (\ref{correlated_cause_specific_Gamma_Model}) is identifiable within the family $\mathcal{H}_1\times\mathcal{H}_2\times\boldsymbol{\Theta}$ if, for some $\boldsymbol{h^{(k)}_{0}},\boldsymbol{\tilde{h}^{(k)}_{0}} \in \mathcal{H}_{k}$, for $k=1,2$, and $\boldsymbol{\theta},\boldsymbol{\tilde{\theta}} \in \boldsymbol{\Theta},$, where $\boldsymbol{\theta}= (\sigma_{1j},\sigma_{2j},\rho_j;\ j=1,\cdots,L),\ \boldsymbol{\tilde{\theta}}= (\tilde{\sigma}_{1j},\tilde{\sigma}_{2j},\tilde{\rho}_j;\ j=1,\cdots,L),$ the equality of the joint sub-distribution functions 
$$F_{j_1,j_2}(t_1,t_2;\boldsymbol{h^{(1)}_{0}},\boldsymbol{h^{(2)}_{0}};\boldsymbol{\theta})= F_{j_1,j_2}(t_1,t_2;\boldsymbol{\tilde{h}^{(1)}_{0}},\boldsymbol{\tilde{h}^{(2)}_{0}};\boldsymbol{\tilde{\theta}}),$$
for all $t_k>0,\ j_k=1,\cdots,L$ and $k=1,2$, implies 
$$h^{(k)}_{0j}(x) = \tilde{h}^{(k)}_{0j}(x) \, \, \text{for all} \, \, x>0,\ k=1,2, \, \, \text{and} \, \, \sigma_{1j} = \tilde{\sigma}_{1j}, \, \sigma_{2j} = \tilde{\sigma}_{2j}, \, \rho_{j} = \tilde{\rho}_{j}, \, \, \text{for} \, \ j= 1,\cdots,L.$$
\end{definition} 

\begin{theorem}\label{Theorem 4.4}
The correlated cause-specific Gamma frailty model (\ref{correlated_cause_specific_Gamma_Model}) with non-parametric baseline cause-specific hazard functions is identifiable within the family $\mathcal{H}_1\times\mathcal{H}_2\times\boldsymbol{\Theta}$, where $\boldsymbol{\Theta}=\big\{\boldsymbol{\theta}=(\sigma_{1j},\sigma_{2j},\rho_j): \sigma_{1j} > 0, \sigma_{2j} > 0,\  0 < \rho_j < \min{(\frac{\sigma_{1j}}{\sigma_{2j}},\frac{\sigma_{2j}}{\sigma_{1j}})},\ j=1,\cdots,L\big\}$, provided $H^{(k)}_{0j}(\cdot)$ is continuous and its inverse, denoted by $H^{-(k)}_{0j}(\cdot)$, exists for all $j=1,\cdots,L$ and $k=1,2$.

\end{theorem}

\begin{proof}

Equality of the unconditional sub-distribution functions $F^{(k)}_{j}(t_{k};\boldsymbol{h^{(k)}_0},\boldsymbol{\theta})$ and $F^{(k)}_{j}(t_{k};\boldsymbol{\tilde{h}^{(k)}_0},\boldsymbol{\tilde{\theta}})$ follows from the condition of equality of the corresponding joint sub-distribution functions in \textbf{Definition} $\mathbf{6}$. Therefore, following similar techniques as those used in the proof of \textbf{Theorem} $\mathbf{6.1}$ for the shared cause-specific Gamma frailty model in the previous section, we get 
\begin{center}
    $h^{(k)}_{0j}(t_k) = \tilde{h}^{(k)}_{0j}(t_k) \, \,\text{and} \, \, \sigma_{j} = \tilde{\sigma}_{j}$
\end{center}
for $t_{k} > 0,\ j = 1,\cdots,L$ and $k = 1,2$. 
In order to prove $\rho_j=\tilde{\rho}_j$, for $j=1,\cdots,L$, we need to consider equality of the unconditional joint sub-distribution functions $F_{j,j}(t_1,t_2;\boldsymbol{h^{(1)}_0},\boldsymbol{h^{(2)}_0},\boldsymbol{\theta})$ and $F_{j,j}(t_1,t_2;\boldsymbol{\tilde{h}^{(1)}_0},\boldsymbol{\tilde{h}^{(2)}_0},\boldsymbol{\tilde{\theta}})$, for all $j = 1,2,\cdots,L$, which also follows from \textbf{Definition} $\mathbf{6}$.\\

Now, the unconditional joint sub-distribution function $F_{j,j}(t_1,t_2;\boldsymbol{h^{(1)}_0},\boldsymbol{h^{(2)}_0},\boldsymbol{\theta})$ can be obtained as 
\begin{align*}
    & \mathbb{E}\Bigg[\int\limits_{0}^{t_1}\int\limits_{0}^{t_2}
    \left(\prod_{K=1}^2 h^{(k)}_{0j}(u_k)\epsilon^{(k)}_{j}\right)
    \times\exp{\bigg(-\sum\limits_{j'=1}^{L} \sum_{k=1}^2
    H^{(k)}_{0j'}(u_k)\epsilon^{(k)}_{j'}\bigg)}du_2du_1\Bigg]\\
    &= \qquad \int\limits_{0}^{t_1}\int\limits_{0}^{t_2}\int\limits_{0}^{\infty}\int\limits_{0}^{\infty}\cdots\int\limits_{0}^{\infty} \left(\prod_{K=1}^2 h^{(k)}_{0j}(u_k)\right)\times 
    \bigg(\frac{\mu^{2}_{0j}}{\mu_{1j}\mu_{2j}}y^{2}_{0j} + \frac{\mu_{0j}}{\mu_{1j}}y_{0j}y_{2j} + \frac{\mu_{0j}}{\mu_{2j}}y_{0j}y_{1j} + y_{1j}y_{2j}\bigg) \\
    &\qquad\times \exp{\bigg[-\sum\limits_{j' (\ne j)=1}^{L} \sum_{k=1}^2 H^{(k)}_{0j'}(u_k)\bigg(\frac{\mu_{0j'}}{\mu_{kj'}}y_{0j'} + y_{kj'}\bigg)\bigg]}\\
    &\qquad\times \exp{\bigg[- \sum_{k=1}^2 
    H^{(k)}_{0j}(u_k)\bigg(\frac{\mu_{0j}}{\mu_{kj}}y_{0j} + y_{kj}\bigg)\bigg]}\\
    &\qquad\times \left(\prod\limits_{l=0}^{2} \prod\limits_{j'=1}^{L} \frac{\mu^{\kappa_{lj'}}_{lj'}}{\Gamma(\kappa_{lj'})}e^{-\mu_{lj'}y_{lj'}} y^{\kappa_{lj'} - 1}_{lj'}dy_{lj'}\right) du_2du_1 \\
    &= \int\limits_{0}^{t_1}\int\limits_{0}^{t_2} \left(\prod_{k=1}^2 h^{(k)}_{0j}(u_k)\right)\times
    \big[I_{1j}(u_1,u_2) + I_{2j}(u_1,u_2) + I_{3j}(u_1,u_2) + I_{4j}(u_1,u_2)\big]du_2du_1, \, \text{(say)},
\end{align*}
where
\begin{align*}
    I_{1j}(u_1,u_2) &= \frac{\mu^{2}_{0j}}{\mu_{1j}\mu_{2j}}\int\limits_{0}^{\infty}\cdots\int\limits_{0}^{\infty} 
    \exp{\bigg[-\sum\limits_{j' (\ne j) =1}^{L}\sum_{k=1}^2
    H^{(k)}_{0j'}(u_k)\bigg(\frac{\mu_{0j'}}{\mu_{kj'}}y_{0j'} + y_{kj'}\bigg)\bigg]}\\
    &\qquad\times\exp{\bigg[-\sum_{k=1}^2 
    H^{(k)}_{0j}(u_k)\bigg(\frac{\mu_{0j}}{\mu_{kj}}y_{0j} + y_{kj}\bigg)\bigg]}y^{2}_{0j} \\
    &\qquad\times\prod\limits_{l=0}^{2}\prod\limits_{j'=1}^{L} \bigg(\frac{\mu^{\kappa_{lj'}}_{lj'}}{\Gamma(\kappa_{lj'})}e^{-\mu_{lj'}y_{lj'}} y^{\kappa_{lj'} - 1}_{lj'}dy_{lj'}\bigg)\\
    &= \frac{\mu^{2}_{0j}}{\mu_{1j}\mu_{2j}}\prod\limits_{j' (\ne j)=1}^{L}\frac{\mu^{\kappa_{0j'}}_{0j'}}{\Gamma(\kappa_{0j'})}\int\limits_{0}^{\infty}\exp{\bigg[-\bigg(1+ \sum_{k=1}^2\frac{H^{(k)}_{0j'}(u_k)}{\mu_{kj'}}\bigg) \mu_{0j'}y_{0j'}\bigg]}y^{\kappa_{0j'} - 1}_{0j'}dy_{0j'}\\
    &\qquad\times \frac{\mu^{\kappa_{0j}}_{0j}}{\Gamma(\kappa_{0j})}\int\limits_{0}^{\infty}\exp{\bigg[-\bigg(1+\sum_{k=1}^2\frac{H^{(k)}_{0j}(u_k)}{\mu_{kj}}\bigg) \mu_{0j}y_{0j}\bigg]}y^{\kappa_{0j} + 1}_{0j}dy_{0j}\\
    &\qquad\times \prod\limits_{j'=1}^{L}\prod_{k=1}^2 \left( \frac{\mu^{\kappa_{kj'}}_{kj'}}{\Gamma(\kappa_{kj'})}\int\limits_{0}^{\infty}\exp{\bigg[-\bigg(1+\frac{H^{(k)}_{0j'}(u_k)}{\mu_{kj'}}\bigg) \mu_{kj'}y_{kj'}\bigg]}y^{\kappa_{kj'} - 1}_{kj'}dy_{kj'} \right)
\end{align*}
\begin{align*}
    &= \frac{\kappa_{0j}(\kappa_{0j} + 1)}{\mu_{1j}\mu_{2j}} 
    \bigg(1+\sum_{k=1}^2\frac{H^{(k)}_{0j}(u_k)}{\mu_{kj}}\bigg)^{-2} \times \\
   &\qquad \prod\limits_{j'=1}^{L}\bigg[\bigg(1+\sum_{k=1}^2 
   \frac{H^{(k)}_{0j'}(u_k)}{\mu_{kj'}}\bigg)^{-\kappa_{0j'}}\times 
   \prod_{k=1}^2\bigg(1+\frac{H^{(k)}_{0j'}(u_k)}{\mu_{kj'}}\bigg)^{-\kappa_{kj'}}\bigg], 
\end{align*}
\begin{align*}
    I_{2j}(u_1,u_2) &= \frac{\mu_{0j}}{\mu_{1j}}\int\limits_{0}^{\infty}\cdots\int\limits_{0}^{\infty}\exp{\bigg[-\sum\limits_{j' (\ne j)=1}^{L}\sum_{k=1}^2
    H^{(k)}_{0j'}(u_k)\bigg(\frac{\mu_{0j'}}{\mu_{kj'}}y_{0j'} + y_{kj'}\bigg)\bigg]}\\
    &\qquad\times \exp{\bigg[-\sum_{k=1}^2
    H^{(k)}_{0j}(u_k)\bigg(\frac{\mu_{0j}}{\mu_{kj}}y_{0j} + y_{kj}\bigg)\bigg]}y_{0j}y_{2j}   \\
    &\qquad\times \bigg( \prod\limits_{l=0}^{2}\prod\limits_{j'=1}^{L} \frac{\mu^{\kappa_{lj'}}_{lj'}}{\Gamma(\kappa_{lj'})}e^{-\mu_{lj'}y_{lj'}} y^{\kappa_{lj'} - 1}_{lj'} dy_{lj'}\bigg) \\
    &= \frac{\mu_{0j}}{\mu_{1j}}\prod\limits_{j' (\ne j)=1}^{L} \frac{\mu^{\kappa_{0j'}}_{0j'}}{\Gamma(\kappa_{0j'})}\int\limits_{0}^{\infty}\exp{\bigg[-\bigg(1+\sum_{k=1}^2 \frac{H^{(k)}_{0j'}(u_k)}{\mu_{kj'}}\bigg) \mu_{0j'}y_{0j'}\bigg]}y^{\kappa_{0j'} - 1}_{0j'}dy_{0j'}\\
    &\qquad\times \frac{\mu^{\kappa_{0j}}_{0j}}{\Gamma(\kappa_{0j})}\int\limits_{0}^{\infty}\exp{\bigg[-\bigg(1+ \sum_{k=1}^2 \frac{H^{(k)}_{0j}(u_k)}{\mu_{kj}}  \bigg)\mu_{0j}y_{0j}\bigg]}y^{\kappa_{0j}}_{0j}dy_{0j}\\
&\qquad\times\prod\limits_{j'=1}^{L}\frac{\mu^{\kappa_{1j'}}_{1j'}}{\Gamma(\kappa_{1j'})}\int\limits_{0}^{\infty}\exp{\bigg[-\bigg(1+\frac{H^{(1)}_{0j'}(u_1)}{\mu_{1j'}}\bigg) \mu_{1j'}y_{1j'} \bigg]}y^{\kappa_{1j'} - 1}_{1j'}dy_{1j'}\\
&\qquad\times \prod\limits_{j' (\ne j)=1}^{L} \frac{\mu^{\kappa_{2j'}}_{2j'}}{\Gamma(\kappa_{2j'})}\int\limits_{0}^{\infty}\exp{\bigg[-\bigg(1+\frac{H^{(2)}_{0j'}(u_2)}{\mu_{2j'}}\bigg) \mu_{2j'}y_{2j'}\bigg]}y^{\kappa_{2j'} - 1}_{2j'}dy_{2j'} \\ 
    &\qquad\times \frac{\mu^{\kappa_{2j}}_{2j}}{\Gamma(\kappa_{2j})}\int\limits_{0}^{\infty}\exp{\bigg[-\bigg(1+\frac{H^{(2)}_{0j}(u_2)}{\mu_{2j}}\bigg) \mu_{2j}y_{2j}\bigg]}y^{\kappa_{2j}}_{2j}dy_{2j} \\
    &= \frac{\kappa_{0j}\kappa_{2j}}{\mu_{1j}\mu_{2j}} 
    \bigg(1+ \sum_{k=1}^2 \frac{H^{(k)}_{0j}(u_k)}{\mu_{kj}} \bigg)^{-1}\times \bigg(1+\frac{H^{(2)}_{0j}(u_2)}{\mu_{2j}} \bigg)^{-1} \\
    &\qquad\times\prod\limits_{j'=1}^{L}\bigg[ \bigg(1+\sum_{k=1}^2\frac{H^{(k)}_{0j'}(u_k)}{\mu_{kj'}}\bigg) ^{-\kappa_{0j'}} \prod_{k=1}^2\bigg(1+\frac{H^{(k)}_{0j'}(u_k)}{\mu_{kj'}}\bigg)^{-\kappa_{kj'}}\bigg], 
\end{align*}
\begin{align*}
    I_{3j}(u_1,u_2) &= \frac{\mu_{0j}}{\mu_{2j}}\int\limits_{0}^{\infty}\cdots\int\limits_{0}^{\infty}\exp{\bigg[-\sum\limits_{j' (\ne j)=1}^{L}\sum_{k=1}^2  H^{(k)}_{0j'}(u_k)\bigg(\frac{\mu_{0j'}}{\mu_{kj'}}y_{0j'} + y_{kj'}\bigg)\bigg]} \\
    &\qquad\times \exp{\bigg[-\sum_{k=1}^2 H^{(k)}_{0j}(u_k)\bigg(\frac{\mu_{0j}}{\mu_{kj}}y_{0j} + y_{kj}\bigg)\bigg]} y_{0j}y_{1j} \\
    &\qquad\times \bigg(\prod\limits_{l=0}^{2} \prod\limits_{j'=1}^{L} \frac{\mu^{\kappa_{lj'}}_{lj'}}{\Gamma(\kappa_{lj'})}e^{-\mu_{lj'}y_{lj'}} y^{\kappa_{lj'} - 1}_{lj'}dy_{lj'}\bigg) \\ 
    &= \frac{\mu_{0j}}{\mu_{2j}}\prod\limits_{j' (\ne j)=1}^{L} \frac{\mu^{\kappa_{0j'}}_{0j'}}{\Gamma(\kappa_{0j'})}\int\limits_{0}^{\infty}\exp{\bigg[-\bigg(1+\sum_{k=1}^2\frac{H^{(k)}_{0j'}(u_k)}{\mu_{kj'}}\bigg) \mu_{0j'}y_{0j'}\bigg]}y^{\kappa_{0j'} - 1}_{0j'}dy_{0j'} \\
    &\qquad\times \frac{\mu^{\kappa_{0j}}_{0j}}{\Gamma(\kappa_{0j})}\int\limits_{0}^{\infty}\exp{\bigg[-\bigg(1+\sum_{k=1}^2\frac{H^{(k)}_{0j}(u_k)}{\mu_{kj}}\bigg) \mu_{0j}y_{0j}\bigg]}y^{\kappa_{0j}}_{0j}dy_{0j} \\
    &\qquad\times \prod\limits_{j' (\ne j)=1}^{L} \frac{\mu^{\kappa_{1j'}}_{1j'}}{\Gamma(\kappa_{1j'})}\int\limits_{0}^{\infty}\exp{\bigg[-\bigg(1+\frac{H^{(1)}_{0j'}(u_1)}{\mu_{1j'}}\bigg) \mu_{1j'}y_{1j'}\bigg]}y^{\kappa_{1j'} - 1}_{1j'}dy_{1j'}  \\
  &\qquad\times \frac{\mu^{\kappa_{1j}}_{1j}}{\Gamma(\kappa_{1j})}\int\limits_{0}^{\infty}\exp{\bigg[-\bigg(1+\frac{H^{(1)}_{0j}(u_1)}{\mu_{1j}}\bigg) \mu_{1j}y_{1j}\bigg]}y^{\kappa_{1j}}_{1j}dy_{1j} \\
&\qquad\times\prod\limits_{j'=1}^{L}\frac{\mu^{\kappa_{2j'}}_{2j'}}{\Gamma(\kappa_{2j'})}\int\limits_{0}^{\infty}\exp{\bigg[-\bigg(1+\frac{H^{(2)}_{0j'}(u_2)}{\mu_{2j'}}\bigg) \mu_{2j'}y_{2j'}\bigg]}y^{\kappa_{2j'} - 1}_{2j'}dy_{2j'} \\
&= \frac{\kappa_{0j}\kappa_{1j}}{\mu_{1j}\mu_{2j}} \bigg(1+ \sum_{k=1}^2 \frac{H^{(k)}_{0j}(u_k)}{\mu_{kj}}\bigg)^{-1}  \bigg(1+\frac{H^{(1)}_{0j}(u_1)}{\mu_{1j}}\bigg)^{-1} \\
&\qquad\times \prod\limits_{j'=1}^{L} \bigg[ \bigg(1+\sum_{k=1}^2\frac{H^{(k)}_{0j'}(u_k)}{\mu_{kj'}}\bigg)^{-\kappa_{0j'}}\prod_{k=1}^2\bigg(1+\frac{H^{(k)}_{0j'}(u_k)}{\mu_{kj'}}\bigg)^{-\kappa_{kj'}}\bigg]
\end{align*}
and 
\begin{align*}
    I_{4j}(u_1,u_2) &=  \int\limits_{0}^{\infty}\cdots\int\limits_{0}^{\infty} \exp{\bigg[-\sum\limits_{j' (\ne j)=1}^{L}\sum_{k=1}^2  H^{(k)}_{0j'}(u_k)\bigg(\frac{\mu_{0j'}}{\mu_{kj'}}y_{0j'} + y_{kj'}\bigg)\bigg]} \\
    &\qquad\times \exp{\bigg[-\sum_{k=1}^2H^{(k)}_{0j}(u_k)\bigg(\frac{\mu_{0j}}{\mu_{kj}}y_{0j} + y_{kj}\bigg)\bigg]} y_{1j}y_{2j} \\
    &\qquad\times \bigg(\prod\limits_{l=0}^{2} \prod\limits_{j'=1}^{L} \frac{\mu^{\kappa_{lj'}}_{lj'}}{\Gamma(\kappa_{lj'})}e^{-\mu_{lj'}y_{lj'}} y^{\kappa_{lj'} - 1}_{lj'}dy_{li}\bigg) 
\end{align*}
\begin{align*}
   &= \prod\limits_{j'=1}^{L}\frac{\mu^{\kappa_{0j'}}_{0j'}}{\Gamma(\kappa_{0j'})}\int\limits_{0}^{\infty}\exp{\bigg[-\bigg(1+\sum_{k=1}^2\frac{H^{(k)}_{0j'}(u_k)}{\mu_{kj'}}\bigg) \mu_{0j'}y_{0j'}\bigg]}y^{\kappa_{0j'} - 1}_{0j'}dy_{0j'} \\
    &\qquad\times \prod\limits_{j' (\ne j)=1}^{L} \frac{\mu^{\kappa_{1j'}}_{1j'}}{\Gamma(\kappa_{1j'})}\int\limits_{0}^{\infty}\exp{\bigg[-\bigg(1+\frac{H^{(1)}_{0j'}(u_1)}{\mu_{1j'}}\bigg) \mu_{1j'}y_{1j'}\bigg]}y^{\kappa_{1j'} - 1}_{1j'}dy_{1j'} \\
    &\qquad\times \frac{\mu^{\kappa_{1j}}_{1j}}{\Gamma(\kappa_{1j})}\int\limits_{0}^{\infty}\exp{\bigg[-\bigg(1+\frac{H^{(1)}_{0j}(u_1)}{\mu_{1j}}\bigg) \mu_{1j}y_{1j}\bigg]}y^{\kappa_{1j}}_{1j}dy_{1j} \\ 
    &\qquad\times \prod\limits_{j' (\ne j)=1}^{L} \frac{\mu^{\kappa_{2j'}}_{2j'}}{\Gamma(\kappa_{2j'})}\int\limits_{0}^{\infty}\exp{\bigg[-\bigg(1+\frac{H^{(2)}_{0j'}(u_2)}{\mu_{2j'}}\bigg) \mu_{2j'}y_{2j'}\bigg]}y^{\kappa_{2j'} - 1}_{2j'}dy_{2j'} \\ 
    &\qquad\times \frac{\mu^{\kappa_{2j}}_{2j}}{\Gamma(\kappa_{2j})}\int\limits_{0}^{\infty}\exp{\bigg[-\bigg(1+\frac{H^{(2)}_{0j}(u_2)}{\mu_{2j}}\bigg) \mu_{2j}y_{2j}\bigg]}y^{\kappa_{2j}}_{2j}dy_{2j} \\
   &= \frac{\kappa_{1j}\kappa_{2j}}{\mu_{1j}\mu_{2j}} \prod_{k=1}^2 
   \bigg(1+\frac{H^{(k)}_{0j}(u_k)}{\mu_{kj}}\bigg)^{-1} \\
   &\qquad\times \prod\limits_{j'=1}^{L}\bigg[ \bigg(1+\sum_{k=1}^2\frac{H^{(k)}_{0j'}(u_k)}{\mu_{kj'}}\bigg)^{-\kappa_{0j'}}\prod_{k=1}^2 \bigg(1+\frac{H^{(k)}_{0j'}(u_k)}{\mu_{kj'}}\bigg)^{-\kappa_{kj'}}\bigg],
\end{align*}
for $j = 1,\cdots,L$. Therefore, the joint sub-distribution function $F_{j,j}(t_1,t_2;\boldsymbol{h^{(1)}_0},\boldsymbol{h^{(2)}_0},\boldsymbol{\theta})$ is given by 
\begin{align*}
    & \int\limits_{0}^{t_1}\int\limits_{0}^{t_2}
    \left(\prod_{k=1}^2\frac{h^{(k)}_{0j}(u_k)}{\mu_{kj}}\right)
    \prod\limits_{j'=1}^{L}\bigg[\bigg(1+\sum_{k=1}^2 \frac{H^{(k)}_{0j'}(u_k)}{\mu_{kj'}}\bigg)^{-\kappa_{0j'}} \prod_{k=1}^2\bigg(1+\frac{H^{(k)}_{0j'}(u_k)}{\mu_{kj'}}\bigg) ^{-\kappa_{kj'}}\bigg] \\
    &\qquad\times \Bigg[\kappa_{0j}(1 + \kappa_{0j})\bigg(1+ \sum_{k=1}^2\frac{H^{(k)}_{0j}(u_k)}{\mu_{kj}}\bigg)^{-2} + \bigg(1+\sum_{k=1}^2\frac{H^{(k)}_{0j}(u_k)}{\mu_{kj}}\bigg)^{-1}\times \sum_{k=1}^2 
    \kappa_{0j}\kappa_{kj}\bigg(1+\frac{H^{(k)}_{0j}(u_k)}{\mu_{kj}}\bigg)^{-1} \\
    &\qquad + \prod_{k=1}^2 \kappa_{kj}\bigg(1+\frac{H^{(k)}_{0j}(u_k)}{\mu_{kj}}\bigg)^{-1}\Bigg]du_2du_1\\
    &= \int\limits_{0}^{t_1}\int\limits_{0}^{t_2}\left(\prod_{k=1}^2 \sigma^{2}_{kj}h^{(k)}_{0j}(u_k)\right)\times     \prod\limits_{j'=1}^{L}\bigg(1+\sum_{k=1}^2\sigma^{2}_{kj'}H^{(k)}_{0j'}(u_k)\bigg)^{-\frac{\rho_{j'}}{\sigma_{1j'}\sigma_{2j'}}} \prod_{k=1}^2 \bigg(1+\sigma^{2}_{kj'}H^{(k)}_{0j'}(u_k)\bigg)^{-\big(\frac{1}{\sigma^{2}_{kj'}} - \frac{\rho_{j'}}{\sigma_{1j'}\sigma_{2j'}}\big)} \\
&\qquad\times \Bigg[\frac{\frac{\rho_{j}}{\sigma_{1j}\sigma_{2j}}\Big(1 + \frac{\rho_{j}}{\sigma_{1j}\sigma_{2j}}\Big)}{\Big(1+\sum_{k=1}^2
\sigma^{2}_{kj}H^{(k)}_{0j}(u_k)\Big)^2}
+ \frac{\frac{\rho_{j}}{\sigma_{1j}\sigma_{2j}}}{\Big(1+\sum_{k=1}^2\sigma^{2}_{kj}H^{(k)}_{0j}(u_k)\Big)} \sum_{k=1}^2 \frac{\frac{1}{\sigma^{2}_{kj}} - \frac{\rho_{j}}{\sigma_{1j}\sigma_{2j}}}{\Big(1+\sigma^{2}_{kj}H^{(k)}_{0j}(u_k)\Big)}  \\
&\quad + \prod_{k=1}^2\left(\frac{\frac{1}{\sigma^{2}_{kj}} - \frac{\rho_{j}}{\sigma_{1j}\sigma_{2j}}}{1+\sigma^{2}_{kj}H^{(k)}_{0j}(u_k)}\right)\Bigg]du_2du_1. 
\end{align*}

Note that the equality of the two joint sub-distribution functions, as mentioned above, also implies the equality of the corresponding joint sub-density functions $f_{j,j}(t_1,t_2;\boldsymbol{h^{(1)}_0},\boldsymbol{h^{(2)}_0},\boldsymbol{\theta})$ and $f_{j,j}(t_1,t_2;\boldsymbol{h^{(1)}_0},\boldsymbol{h^{(2)}_0},\boldsymbol{\tilde{\theta}})$, for $j = 1,2,\cdots,L$, which can be obtained by differentiating the joint sub-distribution functions with respect to $t_1$ and $t_2$. Utilizing the already obtained results 
$\displaystyle{h^{(k)}_{0,j}(t_k) = \tilde{h}^{(k)}_{0,j}(t_k)}$, for all $t_{k}>0$, and $\displaystyle{\sigma_{kj} = \tilde{\sigma}_{kj}}$, for $j = 1,\cdots,L$ and $k = 1,2$, in the equality of the corresponding joint sub-density functions, we get 
\begin{align*}
&\prod\limits_{j'=1}^{L}\Bigg[\bigg(1+\sum_{k=1}^2\sigma^{2}_{kj'}H^{(k)}_{0j'}(t_k)\bigg)^{-\frac{\rho_{j'}}{\sigma_{1j'}\sigma_{2j'}}} \prod_{k=1}^2
   \bigg(1+\sigma^{2}_{kj'}H^{(k)}_{0j'}(t_k)\bigg)^{-\big(\frac{1}{\sigma^{2}_{kj'}} - \frac{\rho_{j'}}{\sigma_{1j'}\sigma_{2j'}}\big)}\Bigg] \\
&\qquad\times\Bigg[\frac{\frac{\rho_{j}}{\sigma_{1j}\sigma_{2j}}\Big(1 + \frac{\rho_{j}}{\sigma_{1j}\sigma_{2j}}\Big)}{\Big(1+\sum_{k=1}^2   \sigma^{2}_{kj}H^{(k)}_{0j}(t_k)\Big)^2} + 
\left(\frac{\frac{\rho_{j}}{\sigma_{1j}\sigma_{2j}}}{1+\sum_{k=1}^2 \sigma^{2}_{kj}H^{(k)}_{0j}(t_k)}\right) \sum_{k=1}^2 \frac{
\frac{1}{\sigma^{2}_{kj}} - \frac{\rho_{j}}{\sigma_{1j}\sigma_{2j}}}
{1+\sigma^{2}_{kj}H^{(k)}_{0j}(t_k)} + \\
&\qquad \prod_{k=1}^2 \left(\frac{\frac{1}{\sigma^{2}_{kj}} - \frac{\rho_{j}}{\sigma_{1j}\sigma_{2j}}}{1+\sigma^{2}_{kj}H^{(k)}_{0j}(t_k)}\right)\Bigg]
\end{align*}
=
\begin{align*}
&\prod\limits_{j'=1}^{L}\Bigg[\bigg(1+\sum_{k=1}^2\sigma^{2}_{kj'}H^{(k)}_{0j'}(t_k)\bigg)^{-\frac{\tilde{\rho}_{j'}}{\sigma_{1j'}\sigma_{2j'}}} \prod_{k=1}^2
   \bigg(1+\sigma^{2}_{kj'}H^{(k)}_{0j'}(t_k)\bigg)^{-\big(\frac{1}{\sigma^{2}_{kj'}} - \frac{\tilde{\rho}_{j'}}{\sigma_{1j'}\sigma_{2j'}}\big)}\Bigg] \\
&\qquad\times\Bigg[\frac{\frac{\tilde{\rho}_{j}}{\sigma_{1j}\sigma_{2j}}\Big(1 + \frac{\tilde{\rho}_{j}}{\sigma_{1j}\sigma_{2j}}\Big)}{\Big(1+\sum_{k=1}^2   \sigma^{2}_{kj}H^{(k)}_{0j}(t_k)\Big)^2} + 
\left(\frac{\frac{\tilde{\rho}_{j}}{\sigma_{1j}\sigma_{2j}}}{1+\sum_{k=1}^2 \sigma^{2}_{kj}H^{(k)}_{0j}(t_k)}\right) \sum_{k=1}^2 \frac{
\frac{1}{\sigma^{2}_{kj}} - \frac{\tilde{\rho}_{j}}{\sigma_{1j}\sigma_{2j}}}
{1+\sigma^{2}_{kj}H^{(k)}_{0j}(t_k)} + \\
&\qquad \prod_{k=1}^2 \left(\frac{\frac{1}{\sigma^{2}_{kj}} - \frac{\tilde{\rho}_{j}}{\sigma_{1j}\sigma_{2j}}}{1+\sigma^{2}_{kj}H^{(k)}_{0j}(t_k)}\right)\Bigg],
\end{align*}
for $t_{k} > 0,\ j = 1,2,\cdots,L$ and $k = 1,2$. 
Letting limit as $t_{1} \to 0+$ and $t_{2} \to 0+$ on both sides of the above equation, we get $\rho_{j} = \tilde{\rho}_{j}$, for $j = 1,2,\cdots,L$. This completes the proof.
\end{proof}

\section{Concluding remarks}

We have considered modeling of bivariate failure time data with competing risks using four different types of frailty. Modeling of  bivariate failure time with competing risks using frailty has been considered before by Bandeen-Roche and Liang (2002),  and Gorfine and Hsu (2011). But model identifiability, which is of utmost importance, has not been studied. In this paper, we consider non-parametric baseline cause-specific hazard functions and Gamma frailty to investigate the identifiability of the corresponding models. We have proved identifiability of the models under very reasonable and general conditions. Upon proving identifiability of the models with non-parametric baseline cause-specific hazard functions and Gamma frailty, identifiability of such models with some parametric baseline cause-specific hazard functions and Gamma frailty is consequently established whenever the class of such parametric baseline cause-specific hazard functions is identifiable within the corresponding associated parameter space. Details of such identifiability results can be found in Ghosh et al. (2024) for some parametric class of  baseline cause-specific hazard functions and Gamma frailty. 

The four different frailty distributions considered in this paper  cover a variety of dependence scenarios. One may still think of more general frailty depending on time as well (See Gorfine and Hsu, 2011). But use of such general frailty distributions not only complicates the derivations making those almost intractable, it also suffers from lack of proper interpretation and justification in practice. 
 
In this paper, we have only restricted our investigation to different types of Gamma frailties because of their popularity. There are several other frailty distributions like Inverse Gaussian, Log-normal, etc.,  which also provide mathematical tractability. It is of interest to consider these different parametric families of frailty distribution and investigate the corresponding identifiability issues. We intend to take up such investigation in the future. 
 
The shared Gamma frailty model and the shared cause-specific Gamma frailty model can be easily generalized for modeling multivariate failure time with more than two individuals. More generally, these two models can be used for clustered failure time data with unequal number of individuals in different clusters. 
While the shared Gamma frailty model can allow different sets of competing risks for different individuals, the shared cause-specific Gamma frailty model requires same set of competing risks for all individuals. Extension of the correlated Gamma frailty model for modeling multivariate failure time requires similar manipulation to induce dependence in the joint distribution of the multivariate frailty $(\epsilon^{(1)},\cdots,\epsilon^{(p)})$, where $p$ denotes the number of individuals by introducing different independent Gamma variables, the $Y_i$'s (See Section 5). Here also the modeling can be used for clustered failure time data with unequal number of individuals possibly with different sets of competing risks. The correlated cause-specific Gamma frailty model can be similarly generalized to multivariate or clustered failure time data with some complicated manipulation to induce dependence, but requiring the same set of competing risks for all the individuals. 

The joint sub-distribution function $F_{j_1j_2}(t_1,t_2;\boldsymbol{h^{(1)}_0},\boldsymbol{h^{(2)}_0},\boldsymbol{\theta})$ is the fundamental theoretical quantity of interest for bivariate failure time data with competing risks. Therefore, it is imperative that the identifiability of any relevant model is  investigated through this quantity, as has been done in this paper. Noting that any arbitrarily and independently censored bivariate failure time data with competing risks can be modeled by means of this joint sub-distribution function, the identifiability results of this paper allows the corresponding models to be applied for such censored data as well. In particular, we can use these identifiable models for bivariate interval censored data with competing risks and also, more particularly, for bivariate current status data with competing risks. 

Even though identifiability is an important aspect of statistical modeling, making inference based on such identifiable models is a different issue. It is apparent from the identifiability results that one requires large enough data to be able to estimate the model parameters. As noted by Li et al. (2001), there may be certain scenarios where the log-likelihood function based on a mathematically identifiable model is nearly flat, specially near the maxima, or has an irregular surface with multiple maxima, in particular, for finite samples. This gives rise to some numerical problems while maximizing the log-likelihood function (See Ghosh et al, 2024).

\section{Refrences}
Bandeen‐Roche, K., and Liang, K. Y. (2002). Modeling multivariate failure time associations in the presence of a competing risk. Biometrika, 89(2), 299-314.\\
Ghosh, B., Dewanji, A., and Das, S. (2024). Parametric Analysis of Bivariate Current Status data with Competing risks using Frailty model. https://doi.org/10.48550/arXiv.2405.05773\\
Dabrowska, D. M. (1988). Kaplan-Meier estimate on the plane. The Annals of Statistics, 1475-1489.\\
Duchateau, L., and Janssen, P. (2008). The frailty model. New York: Springer Verlag.\\
Gorfine, M., and Hsu, L. (2011). Frailty-based competing risks model for multivariate survival data. Biometrics, 67(2), 415-426.\\
Hougaard, P. (1995). Frailty models for survival data. Lifetime data analysis, 1, 255-273.\\
Iachine, I. A. (2004). Identifiability of bivariate frailty models. Preprint, 5.\\
Johnson, N. L., and Kotz, S. (1975). A vector multivariate hazard rate. Journal of Multivariate Analysis, 5(1), 53-66.\\
Li, C. S., Taylor, J. M., and Sy, J. P. (2001). Identifiability of cure models. Statistics and Probability Letters, 54(4), 389-395.\\
Prenen, L., Braekers, R., and Duchateau, L. (2017). Extending the Archimedean copula methodology to model multivariate survival data grouped in clusters of variable size. Journal of the Royal Statistical Society Series B: Statistical Methodology, 79(2), 483-505.\\
Wienke, A. (2010). Frailty models in survival analysis. Chapman and Hall/CRC.\\
Yashin, A. I., Vaupel, J. W., Iachine, I. A. (1995). Correlated individual frailty: an advantageous approach to survival analysis of bivariate data. Mathematical population studies, 5(2), 145-159.
\end{document}